\documentclass[12pt]{amsart}
\usepackage{amsmath}
\usepackage{amsfonts}
\begin{document}
\title[Schur functions for polynomial sequences of binomial type]
{%%The 
Schur type functions associated with \\ polynomial sequences of binomial type}
\author{Minoru ITOH}
\date{}
\address{Department of Mathematics and Computer Science, 
          Faculty of Science,
          Kagoshima University, Kagoshima 890-0065, Japan}
\email{itoh@sci.kagoshima-u.ac.jp }
\thanks{This research was partially supported by 
JSPS Grant-in-Aid for Young Scientists (B) 17740080.}
\maketitle
\begin{abstract}
   We introduce a class of Schur type 
   functions associated with
   polynomial sequences of binomial type.
   This can be regarded as a generalization of
   the ordinary Schur functions and the factorial Schur functions.
   This generalization satisfies some interesting expansion formulas,
   in which there is a curious duality.
   Moreover this class includes examples which are useful to describe the eigenvalues
   of Capelli type central elements of the universal enveloping algebras
   of classical Lie algebras.
\end{abstract}
\theoremstyle{theorem}
   \newtheorem{theorem}{Theorem}[section]
   \newtheorem{proposition}[theorem]{Proposition}%%[section]
   \newtheorem{lemma}[theorem]{Lemma}%%[section]
   \newtheorem{corollary}[theorem]{Corollary}%%[section]
\theoremstyle{remark}
   \newtheorem*{remark}{Remark}
\numberwithin{equation}{section}
%
%%%%%%%%%%%%%%%%%%%%%%%%%%%%%%%%%%%%%%%%%%%%%%%%%%%%%%%%%%%%%%%%%%%%%%%%%%%%%%%%%%
%
\section*{Introduction}
%
%%%%%%%%%%%%%%%%%%%%%%%%%%%%%%%%%%%%%%%%%%%%%%%%%%%%%%%%%%%%%%%%%%%%%%%%%%%%%%%%%
%
In this article,
we introduce a class of Schur type functions associated with polynomial sequences of binomial type.
%%Namely, suggested by the definition of the ordinary Schur function
Namely, inspired by the definition of the ordinary Schur function
$\det(x^{\lambda_i+N-i}_j) / \det(x^{N-i}_j)$,
we consider the following Schur type function: 
$$
   \det(p_{\lambda_i+N-i}(x_j)) / \det(p_{N-i}(x_j)).
$$
Here $\{ p_n(x) \}_{n \geq 0}$ is a polynomial sequence of binomial type.
%%This can be regarded as a generalization of the ordinary Schur functions 
%%and moreover the factorial Schur functions (\cite{BL}, \cite{CL}).
This can be regarded as a generalization of
the ordinary Schur functions and of the factorial Schur functions (\cite{BL}, \cite{CL}).
%%In addition to this, 
%%we also consider the following function:
We also consider the following function:
$$
   \det(p^*_{\lambda_i+N-i}(x_j)) / \det(p^*_{N-i}(x_j)).
$$
Here we put $p^*_n(x) = x^{-1} p_{n+1}(x)$ 
%%(this $p^*_n(x)$ is a polynomial, and satisfies good relations; see Section~\ref{sec:p^*}).
(this is a polynomial, and satisfies good relations; see Section~\ref{sec:p^*}).
The main results of this article 
are some expansion formulas for these functions
and their mysterious duality corresponding to the exchange $p_n(x) \leftrightarrow p^*_n(x)$ 
and the conjugation of partitions
%%(Sections~3, 4, 5, and~6).
(Sections~3--6).
Most of them are proved by elementary and straightforward calculations. 
%%Beside these results, 
Besides these results, 
we also give an application to representation theory of Lie algebras
(Section~8). 

Let us briefly explain this application.
The factorial Schur functions 
are useful to express the eigenvalues of Capelli type central elements 
of the universal enveloping algebras of the general linear Lie algebra
(more precisely, we should say that the ``shifted Schur functions'' are useful;
by the shift of variables,
the factorial Schur functions
are transformed into the shifted Schur functions (\cite{OO1}, \cite{O})).
In this article, we aim to introduce similar Schur type functions
which is useful to express the eigenvalues of Capelli type central elements 
of the universal enveloping algebras of the orthogonal and symplectic Lie algebras.
This aim is achieved in the case of the polynomial sequence 
corresponding to the central difference.
This case associated with the central difference 
is also related with 
the analogues of the shifted Schur functions given in \cite{OO2},
which were introduced with a similar aim.
Moreover, from this investigation of eigenvalues,
we see the relation
between following two central elements of $U(\mathfrak{o}_N)$ and $U(\mathfrak{sp}_N)$
(Theorem~\ref{thm:relation_with_analogues_of_quantum_immanants}):
(a) the central elements in terms of the column-determinant and 
the column-permanent given in \cite{W} and \cite{I6},
(b) analogues of the quantum immanants given in \cite{OO2}.
Our class is a natural generalization of the Schur functions
containing these interesting functions related with
the classical Lie algebras.

Various generalizations are known for the Schur functions.
Many of them are obtained by replacing the ordinary powers 
by some polynomial sequence
(further generalizations are known; see \cite{M2}).
In particular, the generalization associated with the polynomials 
in the form $p_n(x) = \prod_{k=1}^n (x-a_k)$ is well known \cite{M1},
and this contains the factorial Schur function.
In this article, we consider another generalization
which is not particularly large but includes interesting phenomena and examples.
%
%%%%%%%%%%%%%%%%%%%%%%%%%%%%%%%%%%%%%%%%%%%%%%%%%%%%%%%%%%%%%%%%%%%%%%%%%%%%%%%%%%%
%
\section{Polynomial sequences of binomial type}
%
%%%%%%%%%%%%%%%%%%%%%%%%%%%%%%%%%%%%%%%%%%%%%%%%%%%%%%%%%%%%%%%%%%%%%%%%%%%%%%%%%%%
%
First, we recall the properties of polynomial sequences of binomial type.
See \cite{MR}, \cite{R}, \cite{RKO}, and \cite{S} for further details.
\subsection{{}}
We start with the definition.
%%A polynomial sequence $\{ p_n(x) \}_{n \geq 0}$ 
%%in which the degree of each polynomial is equal to its index,
%%is said to be of binomial type,
%%when the following relation holds for any $n \geq 0$:
A polynomial sequence $\{ p_n(x) \}_{n \geq 0}$ 
in which the degree of each polynomial is equal to its index
is said to be {\it of binomial type}
when the following relation holds for any $n \geq 0$:
$$
   p_n(x+y) = \sum_{k \geq 0} {n \choose k} p_k(x) p_{n-k}(y).
$$

Let us see some examples.
First, the sequence $\{ x^n \}_{n \geq 0}$ of the ordinary powers
is of binomial type,
because we have the relation 
$$
   (x+y)^n = \sum_{k \geq 0} {n \choose k} x^k y^{n-k}
$$
(the ordinary binomial expansion).
As other typical examples,
some factorial powers are well known.
We define the rising factorial power $x^{\overline{n}}$ 
and falling factorial power $x^{\underline{n}}$ by 
$$
   x^{\overline{n}} = x (x+1) \cdots (x+n-1), \qquad
   x^{\underline{n}} = x (x-1) \cdots (x-n+1).
$$
Then $\{ x^{\overline{n}} \}_{n \geq 0}$ and $\{ x^{\underline{n}} \}_{n \geq 0}$ are
also of binomial type. 
Indeed the following relations hold (\cite{MR}, \cite{RKO}):
$$
   (x+y)^{\overline{n}} 
   = \sum_{k \geq 0} {n \choose k} x^{\overline{k}} y^{\overline{n-k}}, \qquad
   (x+y)^{\underline{n}} 
   = \sum_{k \geq 0} {n \choose k} x^{\underline{k}} y^{\underline{n-k}}.
$$

It is easily seen that $p_n(0) = \delta_{n,0}$, 
when  $\{ p_n(x) \}_{n \geq 0}$ is of binomial type.

%%%%%%%%%%%%%%%%%%%%%%%%%%%%%%%%%%%%%%%%%%%%%%%%%%%%%%%%%%%%%%%%%%%%%%%%%%%%%%
%
\subsection{{}}
A natural correspondence is known between
polynomial sequences of binomial type and
delta operators \cite{RKO}.

We recall the definition of delta operators.
A linear operator $Q = Q_x \colon \mathbb{C}[x] \to \mathbb{C}[x]$
%%is called a ``delta operator,''
is called a {\it delta operator}
when the following two properties hold:
(i) $Q$ reduces degrees of polynomials by one;
(ii) $Q$ is shift-invariant
%%(namely, $Q$ commutes with
%%all shift operators $E^a \colon f(x) \mapsto f(x+a)$).
(that is, $Q$ commutes with
all shift operators $E^a \colon f(x) \mapsto f(x+a)$).
A typical example is 
the differentiation $D=\frac{d}{dx}$.
%%Moreover the forward difference $\Delta^+$
%%and the backward difference $\Delta^-$ are delta operators:
Moreover the {\it forward difference} $\Delta^+$
and the {\it backward difference} $\Delta^-$ are delta operators:
$$
   \Delta^+ \colon f(x) \mapsto f(x+1) - f(x),\qquad
   \Delta^- \colon f(x) \mapsto f(x) - f(x-1).
$$
Every delta operator can be written 
%%as a power series of the differentiation $D$ in the following form
as a power series in the differentiation operator $D$ of the following form
with $a_1,a_2,\ldots \in \mathbb{C}$, $a_1 \ne 0$:
$$
   Q = a_1 D + a_2 D^2 + a_3 D^3 + \cdots.
$$

There is a natural one-to-one correspondence between these delta operators and polynomial sequences of binomial type.
These are related via the relation
\begin{equation}
   Q p_n(x) = n p_{n-1}(x). 
\label{eq:relation_between_op_and_seq}
\end{equation}
Namely, for a polynomial sequence  of binomial type $\{ p_n(x) \}_{n \geq 0}$, 
the linear operator $Q \colon \mathbb{C}[x] \to \mathbb{C}[x]$
determined by (\ref{eq:relation_between_op_and_seq}) is a delta operator.
Conversely, for any delta operator $Q$, 
a polynomial sequence $\{ p_n(x) \}_{n \geq 0}$ is uniquely determined by (\ref{eq:relation_between_op_and_seq}) 
and the relation $p_n(0) = \delta_{n,0}$, 
%%and this $\{ p_n(x) \}_{n \geq 0}$ is of binomial type
and this sequence is of binomial type
(these are called {\it basic polynomials}).

For example, the differentiation $D =\frac{d}{dx}$ 
corresponds to the sequence $\{ x^n \}_{n \geq 0}$, 
because $D x^n = n x^{n-1}$.
Similarly, the forward difference $\Delta^+$
and the backward difference $\Delta^-$
correspond to the sequences $\{ x^{\underline{n}} \}_{n \geq 0}$ 
and $\{ x^{\overline{n}} \}_{n \geq 0}$, respectively:
$$
   \Delta^- x^{\overline{n}} = n x^{\overline{n-1}}, \qquad
   \Delta^+ x^{\underline{n}} = n x^{\underline{n-1}}.
$$

%%%%%%%%%%%%%%%%%%%%%%%%%%%%%%%%%%%%%%%%%%%%%%%%%%%%%%%%%%%%%%%%%%%%%%%%%%%%%%%
%
\subsection{{}}\label{sec:p^*}
There is another interesting polynomial sequence associated with a polynomial sequence of binomial type. 
For a polynomial sequence  $\{ p_n(x) \}_{n \geq 0}$ of binomial type,
we put
$$
   p^*_n(x) = x^{-1} p_{n+1}(x).
$$
This is a polynomial, because the constant term of $p_{n+1}(x)$ 
is equal to $0$ if $n \geq 0$ as seen above.
This $p^*_n(x)$ satisfies the following relations (we can prove this by induction).
In other words, $\{ p^*_n(x) \}_{n \geq 0}$ is a Sheffer sequence \cite{R}.

\begin{proposition} \sl We have
   $$
      Q p^*_n(x)  = n p^*_{n-1}(x), \qquad
	  p^*_n(x+y)  = \sum_{k \geq 0} {n \choose k} p_k(x) p^*_{n-k}(y) 
	  = \sum_{k \geq 0} {n \choose k} p^*_k(x) p_{n-k}(y).
   $$
\end{proposition}

These polynomials can be extended naturally for $n < 0$
as elements of 
$\mathbb{C}((x^{-1})) 
= \{ \sum_{k \leq n} a_k x^k \,|\, a_k \in \mathbb{C}, \, n \in \mathbb{Z} \}$.
Namely 
we have the following proposition
(this is also proved by induction):

\begin{proposition} \sl 
   Let $Q$ be a delta operator.
%%   Then, there uniquely exist $\{ p_n(x) \}_{n \in \mathbb{Z}}$ 
   Then there uniquely exist $\{ p_n(x) \}_{n \in \mathbb{Z}}$ 
   and $\{ p^*_n(x) \}_{n \in \mathbb{Z}}$
   satisfying the relations
   $$
      Q p_n(x) = n p_{n-1}(x), \qquad
	  Q p^*_n(x) = n p^*_{n-1}(x), \qquad
	  p_{n+1}(x) = x p^*_n(x)
   $$
   and
   \begin{align*}
	  p_n(x+y) & = \sum_{k \geq 0} {n \choose k} p_k(x) p_{n-k}(y), \\
	  p^*_n(x+y) & = \sum_{k \geq 0} {n \choose k} p_k(x) p^*_{n-k}(y) 
	  = \sum_{k \geq 0} {n \choose k} p^*_k(x) p_{n-k}(y).
   \end{align*}
   Here we regard the last two relations as equalities in $\mathbb{C}[x]((y^{-1}))$.
\end{proposition}

Note that $p^*_{-1}(x)$ must be equal to $x^{-1}$,
because $x p^*_{-1}(x) = p_0(x) = 1$.
Thus, this extension is unique.

From now on, we denote these polynomials associated with the delta operator $Q$
by $p_n(x) = p^Q_n(x)$ and $p^*_n(x) = p^{*Q}_n(x)$.

The polynomial $p^*_n(x)$ is not a mere supplementary object, 
%%but plays as an important role as $p_n$.
but plays a role as important as $p_n$.
%%We will see some dualities between these two polynomials in this article.
We will exhibit some dualities between these two polynomials.
%%%%%%%%%%%%%%%%%%%%%%%%%%%%%%%%%%%%%%%%%%%%%%%%%%%%%%%%%%%%%%%%%%
%
\subsection{{}}
%
%%It is interesting to consider the following operator:
Consider the following operator:
$$
   R_x = [Q_x,x] = Q_x x - x Q_x.
$$
We can easily see that $R_x p^*_k(x) = p_k(x)$,
and $R_x$ is invertible and shift-invariant.
Let us put $p_k^{(a)}(x) = R_x^{a} p_k(x)$ for $a \in \mathbb{Z}$,
so that $p_k^{(-1)}(x) = p^*_k(x)$.
As seen by induction, this $p_k^{(a)}(x)$ satisfies the relation
$$
   p_n^{(a+b)}(x+y) 
   = \sum_{k \geq 0} {n \choose k} p_{n-k}^{(a)}(x) p_{k}^{(b)}(y).
$$
In the case of $Q = \Delta^+$, this operator $R_x$ maps $f(x)$ to $f(x+1)$
%%(namely $p^{(a)}_n(x) = (x+a-1)^{\underline{n}}$),
(that is, $p^{(a)}_n(x) = (x+a-1)^{\underline{n}}$),
and this is an algebra automorphism on $\mathbb{C}[x]$,
but this is not true for general $Q$.
See also Remark in Section~\ref{sec:definition}.

%%%%%%%%%%%%%%%%%%%%%%%%%%%%%%%%%%%%%%%%%%%%%%%%%%%%%%%%%%%%%%%%%%
%
\subsection{{}}
In the remainder of this article,
we assume that $Q$ is {\it normalized}
in the sense that the coefficient of $D$ is equal to $1$:
$$
   Q = D + a_2 D^2 + a_3 D^3 + \cdots.
$$
Under this assumption, 
the associated polynomials $p_n(x) = p^Q_n(x)$ 
and $p^*_n(x) = p^{*Q}_n(x)$ become monic.
Conversely, the delta operator associated with a monic polynomial sequence of binomial type
is automatically normalized.
Thus the following assumptions are equivalent:
(i)~$Q$ is normalized;
(ii)~$p^Q_n(x)$ is monic;
(iii)~$p^{*Q}_n(x)$ is monic.
This assumption is not an essential one, but merely for simplicity.
%%%%%%%%%%%%%%%%%%%%%%%%%%%%%%%%%%%%%%%%%%%%%%%%%%%%%%%%%%%%%%%%%%
%
\subsection{{}}
Let us see some fundamental examples.

\medskip

\noindent
(1)
In the case $Q = D = \frac{d}{dx}$, we have $p_n(x) = x^n$ and $p^*_n(x) = x^n$. 

\medskip

\noindent
(2)
In the case $Q = \Delta^+$, we have $p_n(x) = x^{\underline{n}}$,
and hence $p^*_n(x) = (x-1)^{\underline{n}}$. 

\medskip

\noindent
(3)
In the case $Q = \Delta^-$, we have $p_n(x) = x^{\overline{n}}$,
and hence $p^*_n(x) = (x+1)^{\overline{n}}$. 

\medskip

\noindent
(4) 
%%The case of central difference.
We define the {\it central difference} $\Delta^0$ by
$$
   \Delta^0 f(x) = f(x+\tfrac{1}{2}) - f(x-\tfrac{1}{2}).
$$
This is also a delta operator.
In the case $Q = \Delta^0$, we have $p^*_n(x) = x^{\overline{\underline{n}}}$.
Here we put
$$
   x^{\overline{\underline{n}}} 
   = (x + \tfrac{n-1}{2}) (x + \tfrac{n-3}{2}) \cdots (x - \tfrac{n-1}{2}).
$$
Hence, $p_n(x)$ is expressed as $p_n(x) = x \cdot x^{\overline{\underline{n}}}$. 
This is seen by a direct calculation.

\medskip

See \cite{R}, \cite{RKO}, and \cite{S} for other examples 
%%(the Abel polynomials, the Laguerre polynomials, etc.).
(Abel polynomials, Laguerre polynomials, etc.).
%
%%%%%%%%%%%%%%%%%%%%%%%%%%%%%%%%%%%%%%%%%%%%%%%%%%%%%%%%%%%%%%%%%%%%%%%%%%
%
\section{Definition of Schur type functions}\label{sec:definition}
%
%%%%%%%%%%%%%%%%%%%%%%%%%%%%%%%%%%%%%%%%%%%%%%%%%%%%%%%%%%%%%%%%%%%%%%%%%%
%
Let us define our main objects, the Schur type functions associated with
a polynomial sequence of binomial type.
Let $p_n(x) = p^Q_n(x)$ and $p^*_n(x) = p^{*Q}_n(x)$ be 
polynomials corresponding to a normalized delta operator $Q$.

For $\lambda = (\lambda_1,\ldots,\lambda_N) \in \mathbb{Z}^N$, we consider the following determinants:
\begin{align*}
   \tilde{s}^Q_{\lambda}(x_1,\ldots,x_N) & = \det 
   \begin{pmatrix}
       p_{\lambda_1 + N-1}(x_1) & \hdots & p_{\lambda_1 + N-1}(x_N) \\
       p_{\lambda_2 + N-2}(x_1) & \hdots & p_{\lambda_2 + N-2}(x_N) \\
       \vdots & & \vdots \\
	   p_{\lambda_N + 0}(x_1) & \hdots & p_{\lambda_N + 0}(x_N) 
   \end{pmatrix}, \\
   \tilde{s}^{*Q}_{\lambda}(x_1,\ldots,x_N) & = \det 
   \begin{pmatrix}
       p^*_{\lambda_1 + N-1}(x_1) & \hdots & p^*_{\lambda_1 + N-1}(x_N) \\
       p^*_{\lambda_2 + N-2}(x_1) & \hdots & p^*_{\lambda_2 + N-2}(x_N) \\
       \vdots & & \vdots \\
	   p^*_{\lambda_N + 0}(x_1) & \hdots & p^*_{\lambda_N + 0}(x_N) 
   \end{pmatrix}.
\end{align*}
We regard these as elements of
$$
   \mathbb{C}((x_1^{-1},\ldots,x_N^{-1})) = 
   \left\{ \left. 
   \sum_{k_1 \leq n_1,\ldots,k_N \leq n_N} a_{k_1,\ldots,k_N} x_1^{k_1} \cdots x_N^{k_N} 
   \, \right| \, 
   a_{k_1,\ldots,k_N} \in \mathbb{C},\,
   n_1,\ldots,n_N \in \mathbb{Z} 
   \right\}.
$$
%%It is easy to see that these are alternating in $x_1,\ldots,x_N$.
It is easy to see that these two determinants are alternating in $x_1,\ldots,x_N$.
Moreover, 
%%when $\lambda = \emptyset$,
when $\lambda = \emptyset = (0,\ldots,0)$, 
these are equal to the difference product:
\begin{equation}
   \tilde{s}^Q_{\emptyset}(x_1,\ldots,x_N) 
   = 
   \tilde{s}^{*Q}_{\emptyset}(x_1,\ldots,x_N) 
   = \Delta(x_1,\ldots,x_N) = \prod_{1 \leq i<j \leq N} (x_i - x_j).
\label{eq:difference_product}
\end{equation}
%%Indeed we can transform these to the ordinary Vandermonde determinant
Indeed, we can transform these to the ordinary Vandermonde determinant
by elementary row operations,
because $p_n(x)$ and $p^*_n(x)$ are monic.
%%Noting this, we consider the following functions:
Having noted this, we consider the following functions:
\begin{align*}
     s^Q_{\lambda}(x_1,\ldots,x_N) 
     & = \tilde{s}^Q_{\lambda}(x_1,\ldots,x_N) 
   / \tilde{s}^Q_{\emptyset}(x_1,\ldots,x_N), \\
     s^{*Q}_{\lambda}(x_1,\ldots,x_N) 
     & = \tilde{s}^{*Q}_{\lambda}(x_1,\ldots,x_N) 
   / \tilde{s}^{*Q}_{\emptyset}(x_1,\ldots,x_N).
\end{align*} 
If $Q = D$, these are equal to the ordinary Schur functions.
%%We can easily see that these are symmetric functions for any delta operator $Q$,
We can easily see that they are symmetric functions for any delta operator $Q$,
and their highest degree parts are equal to the ordinary Schur function.

From now on, we omit the superscript $Q$, when it is clear from the context.

%%These $s_{\lambda}$ and $s^*_{\lambda}$ are polynomials, 
The functions $s_{\lambda}$ and $s^*_{\lambda}$ are polynomials, 
if $\lambda_1,\ldots,\lambda_N \geq 0$.
When $\lambda = (\lambda_1,\ldots,\lambda_N)$ is a partition, 
namely when $\lambda_1 \geq \cdots \geq \lambda_N \geq 0$, 
we can regard $\lambda$ as a Young diagram.
In this case, the polynomials $s_{\lambda}$ (respectively, $s^*_{\lambda}$)
form a basis of the space of symmetric polynomials.

We can also consider the counterparts of elementary symmetric functions
and complete homogeneous symmetric functions
(note that $h_k$ and $h^*_k$ can be defined even if $k$ is a negative integer):
\begin{align*}
   e_k(x_1,\ldots,x_N) &= s_{(1^k)}(x_1,\ldots,x_N), \\
   e^*_k(x_1,\ldots,x_N) &= s^*_{(1^k)}(x_1,\ldots,x_N), \\
   h_k(x_1,\ldots,x_N) &= s_{(k)}(x_1,\ldots,x_N), \\
   h^*_k(x_1,\ldots,x_N) &= s^*_{(k)}(x_1,\ldots,x_N).
\end{align*}
Here we used the abbreviation
$$
   (a_1^{m_1},\ldots,a_n^{m_n}) 
   = (\overbrace{a_1,\ldots,a_1}^{\text{$m_1$ times}},
   \ldots,\overbrace{a_n,\ldots,a_n}^{\text{$m_n$ times}},0,\ldots,0).
$$
%%These are not equal to the ordinary elementary symmetric functions
These functions are not equal to the ordinary elementary symmetric functions
and the ordinary complete symmetric functions in general,
but there are two exceptions.
Namely $e_N$ and $h^*_{-N}$ are independent of $\{ p_n(x) \}$.
\begin{proposition}\label{prop:special_e_and_h*}\sl
   We have
   $$
      e_N(x_1,\ldots,x_N) = x_1 \cdots x_N, \qquad
	  h^*_{-N}(x_1,\ldots,x_N) = \frac{1}{x_1 \cdots x_N}.
   $$
\end{proposition}
This is easy from the following more general relation:
\begin{proposition} \sl
   We have
   $$
      s_{(\lambda_1,\ldots,\lambda_N)}(x_1,\ldots,x_N) 
	  = s^*_{(\lambda_1 - 1,\ldots,\lambda_N - 1)}(x_1,\ldots,x_N)
	  \cdot x_1 \cdots x_N.
   $$
\end{proposition}
This is immediate from
$$
      \tilde{s}_{(\lambda_1,\ldots,\lambda_N)}(x_1,\ldots,x_N) 
	  = \tilde{s}^*_{(\lambda_1 - 1,\ldots,\lambda_N - 1)}(x_1,\ldots,x_N)
	  \cdot x_1 \cdots x_N.
$$

The following relation between $s$ and $s^*$ is also confirmed by a direct calculation:
\begin{proposition}\label{prop:substitute_0} \sl
   When $\lambda_i \geq 0$, we have
   $$
      s^*_{(\lambda_1,\ldots,\lambda_N)}(x_1,\ldots,x_N) 
	  = s_{(\lambda_1,\ldots,\lambda_N)}(x_1,\ldots,x_N,0).
   $$
\end{proposition}
\begin{remark}
The shifted Schur function is defined as follows \cite{OO1}:
$$
   \det((x_j + N - j)^{\underline{\lambda_i+N-i}}) / 
   \det((x_j + N - j)^{\underline{N-i}}).
$$
In the case of $Q_x = \Delta^+_x$, we have $p_n(x) = x^{\underline{n}}$,
and $R_x = [Q_x, x] = Q_x x - x Q_x$ is equal to the algebra automorphism $f(x) \mapsto f(x+1)$.
Thus this function can be expressed as
\begin{align*}
   & \det(p^{(N-j)}_{\lambda_i+N-i}(x_j)) / \det(p^{(N-j)}_{N-i}(x_j)) \\
   & \qquad
   =
   R_{x_1}^{N-1} R_{x_2}^{N-2} \cdots R_{x_N}^{0} 
   \det(p_{\lambda_i+N-i}(x_j)) / \det(p_{N-i}(x_j)).
\end{align*}
Noting this and Proposition~\ref{prop:substitute_0}, 
we can naturally consider the projective limit of this function.
This is an advantage of considering the shift $x_j \mapsto x_j +N-j$.

How about the case of an arbitrary delta operator $Q_x$?
In general, we do not have such a good relation,
because $R_x$ is not an algebra automorphism
(thus this is not a polynomial in general, even if $\lambda$ is a partition). 
%%Thus it seems not easy to consider a natural infinite-variable version for general $Q_x$.
Thus it does not seem easy to consider a natural infinite-variable version for general $Q_x$.
\end{remark}
%
%%%%%%%%%%%%%%%%%%%%%%%%%%%%%%%%%%%%%%%%%%%%%%%%%%%%%%%%%%%%%%%%%%%%%%%%%%
%
\section{Expansions of Schur type functions}
%
%%%%%%%%%%%%%%%%%%%%%%%%%%%%%%%%%%%%%%%%%%%%%%%%%%%%%%%%%%%%%%%%%%%%%%%%%%
%
For the Schur type functions defined in the previous section,
we have the following expansions 
(this can be regarded as a generalization of Example 10 in Section~I.3 in \cite{M1}):
\begin{theorem}\label{thm:expansion_of_s}
   \sl For $\lambda_1 \geq \lambda_2 \geq \cdots \geq \lambda_N$, we have
   \begin{align*}
      s_{\lambda}(x_1 + u,\ldots,x_N +u)
	  & = \sum_{\mu \subset \lambda} d_{\lambda\mu}(u) s_{\mu}(x_1,\ldots,x_N),\\
      s^*_{\lambda}(x_1 + u,\ldots,x_N +u)
	  & = \sum_{\mu \subset \lambda} d_{\lambda\mu}(u) s^*_{\mu}(x_1,\ldots,x_N),\\
      s^*_{\lambda}(x_1 + u,\ldots,x_N +u)
	  & = \sum_{\mu \subset \lambda} d^*_{\lambda\mu}(u) s_{\mu}(x_1,\ldots,x_N)
   \end{align*}
   as equalities in $\mathbb{C}[u]((x_1^{-1},\ldots,x_N^{-1}))$.
   Here $\mu$ runs over $\mu = (\mu_1,\ldots,\mu_N)$ such that
   $$
      \mu_1 \geq \mu_2 \geq \cdots \geq \mu_N, \qquad 
	  \mu_1 \leq \lambda_1,\ldots,
	  \mu_N \leq \lambda_N,
   $$ 
   and $d_{\lambda\mu}(u)$ and $d^*_{\lambda\mu}(u)$ are defined by
   \begin{align*}
      d_{\lambda\mu}(u) & = \det 
	  \left( 
	     {\lambda_i + N-i \choose \lambda_i - \mu_j - i + j} p_{\lambda_i - \mu_j -i +j} (u) 
	  \right)_{1 \leq i,j \leq N}, \\
      d^*_{\lambda\mu}(u) & = \det 
	  \left( 
	     {\lambda_i + N-i \choose \lambda_i - \mu_j - i + j} p^*_{\lambda_i - \mu_j -i +j} (u) 
	  \right)_{1 \leq i,j \leq N}.
   \end{align*}
\end{theorem}

To prove this, we use the Cauchy--Binet formula:

\begin{proposition}\label{prop:Cauchy--Binet}
   \sl We have
   $$
      \det (AB)_{(i_1,\ldots,i_k),(j_1,\ldots,j_k)} 
	  = \sum_{1 \leq r_1 < \cdots < r_k} 
	  \det A_{(i_1,\ldots,i_k),(r_1,\ldots,r_k)} 
	  \det B_{(r_1,\ldots,r_k),(j_1,\ldots,j_k)}.
   $$
   Here we put $X_{(i_1,\ldots,i_k),(j_1,\ldots,j_k)} = (x_{i_a,j_b})_{1 \leq a,b \leq k}$
   for a matrix $X = (x_{ij})$.
   Note that this holds when the multiplication is defined, 
   even if the sizes of these matrices are infinite.
\end{proposition}

\noindent
{\it Proof of Theorem~{\sl \ref{thm:expansion_of_s}}. }
We put $l_i = \lambda_i + N - i$,
and consider the matrix
$$
   A = \begin{pmatrix}
   p_{l_1} (x_1 + u) & \hdots & p_{l_1} (x_N + u) \\
   \vdots & & \vdots \\
   p_{l_N} (x_1 + u) & \hdots & p_{l_N} (x_N + u) \\   
   \end{pmatrix}.
$$
The $(i,j)$th entry $p_{l_i}(x_j + u)$ can be expanded as
$$
   p_{l_i}(x_j + u)
   = \sum_{k \leq l_1} {l_i \choose l_i - k} p_{l_i - k}(u) p_k(x_j).
$$
Thus $A$ is expressed as $A = BC$
with the $N \times \infty$ matrix $B$ and 
the $\infty \times N$ matrix $C$ defined by
\begin{align*}
   B &= \begin{pmatrix}
   {l_1 \choose 0} p_0(u) & {l_1 \choose 1} p_1(u) & \hdots \\
   {l_2 \choose l_2 - l_1} p_{l_2 - l_1}(u) & {l_2 \choose l_2 - l_1 +1} p_{l_2 - l_1 +1}(u) & \hdots \\
   \vdots & \vdots & \\
   {l_N \choose l_N - l_1} p_{l_N - l_1}(u) & {l_N \choose l_N - l_1 +1} p_{l_N - l_1 +1}(u) & \hdots
   \end{pmatrix},\\
   C &= \begin{pmatrix}
   p_{l_1}(x_1) &  p_{l_1}(x_2) & \hdots & p_{l_1}(x_N) \\
   p_{l_1 - 1}(x_1) &  p_{l_1 - 1}(x_2) & \hdots & p_{l_1 - 1}(x_N) \\
   \vdots & \vdots & & \vdots
   \end{pmatrix}.
\end{align*}
Applying the Cauchy--Binet formula (Proposition~\ref{prop:Cauchy--Binet}) to this, we have
\begin{align*}
   & \tilde{s}_{\lambda}(x_1 + u,\ldots,x_N +u) \\
   & \qquad
   = \det A \\
   & \qquad
   = \sum_{k_1,\ldots,k_N} 
   \det 
   \begin{pmatrix}
   {l_1 \choose l_1 - k_1} p_{l_1 - k_1}(u) & {l_1 \choose l_1 - k_2} p_{l_2 - k_2}(u) & \hdots & {l_1 \choose l_1 - k_N} p_{l_1 - k_N}(u)\\
   {l_2 \choose l_2 - k_1} p_{l_2 - k_1}(u) & {l_2 \choose l_2 - k_2} p_{l_2 - k_2}(u) & \hdots & {l_2 \choose l_2 - k_N} p_{l_2 - k_N}(u)\\
   \vdots & \vdots & & \vdots \\
   {l_N \choose l_N - k_1} p_{l_N - k_1}(u) & {l_N \choose l_N - k_2} p_{l_N - k_2}(u) & \hdots & {l_N \choose l_N - k_N} p_{l_N - k_N}(u)\\
   \end{pmatrix}\\
   & \qquad\qquad\qquad
   \cdot 
   \det 
   \begin{pmatrix}
   p_{k_1}(x_1) &  p_{k_1}(x_2) & \hdots & p_{k_1}(x_N) \\
   p_{k_2}(x_1) &  p_{k_2}(x_2) & \hdots & p_{k_2}(x_N) \\
   \vdots & \vdots & & \vdots \\
   p_{k_N}(x_1) &  p_{k_N}(x_2) & \hdots & p_{k_N}(x_N) 
   \end{pmatrix} \\
   & \qquad
   = \sum_{\mu \subset \lambda} d_{\lambda\mu}(u) \tilde{s}_{\mu}(x_1,\ldots,x_N).
\end{align*}
Here,
%%the first sum runs over $k_1,\ldots,k_N$ satisfying 
the first sum is over $k_1,\ldots,k_N$ satisfying 
$k_i \leq l_i$ and $k_1 > \cdots > k_N$,
%%and the second sum runs over $\mu_1,\ldots,\mu_N$ 
and the second over $\mu_1,\ldots,\mu_N$ 
satisfying $\mu_i \leq \lambda_i$ and $\mu_1 \geq \cdots \geq \mu_N$
%%(namely we define $\mu_i$ by $k_i = \mu_i + N - i$).
(we define $\mu_i$ by $k_i = \mu_i + N - i$).
Dividing this by $\Delta(x_1 + u,\ldots,x_N +u) = \Delta(x_1,\ldots,x_N)$,
we have the assertion.
\hfil\qed

\medskip

%%It is interesting to consider the following variant of this $d_{\lambda\mu}(u)$,
It is interesting to consider the following variant of $d_{\lambda\mu}(u)$
when $\lambda$ and $\mu$ are partitions:
$$
   \hat{d}_{\lambda\mu}(u) =
   \frac{\prod_j (\mu_j + N - j)!}{\prod_i (\lambda_i + N - i)!} d_{\lambda\mu}(u) 
   = \det 
     \left( 
	    p_{(\lambda_i - \mu_j -i +j)} (u) 
	 \right)_{1 \leq i,j \leq N}.
$$
Here we put $p_{(n)}(x) = \frac{1}{n!} p_n(x)$ 
suggested by the notation of divided power $x^{(n)} = \frac{1}{n!}x^n$.
It is easily seen that $\hat{d}_{\lambda\mu}(u)$ is independent of~$N$,
though $d_{\lambda\mu}(u)$ depends on~$N$.
Namely, $\hat{d}_{\lambda\mu}$
does not change, even if we append some zeros at the ends of~$\lambda$ and~$\mu$. 
For this $\hat{d}_{\lambda\mu}(u)$, 
the following duality holds:
\begin{theorem}\label{thm:duality_coeff}\sl
   For two partitions $\lambda$ and $\mu$, 
   we have
   $$
      \hat{d}_{\lambda\mu}(u) 
      = (-)^{|\lambda|-|\mu|}\hat{d}_{\lambda' \mu'}(-u).
   $$
   Here $\lambda'$ and $\mu'$ mean the conjugates of $\lambda$ and $\mu$, respectively.
\end{theorem}
This theorem follows from Theorem~\ref{thm:Cauchy_1} below.
%%Moreover this theorem can be also deduced from the following relation 
%%(this is essentially the same as the relation~(2.9) in~\cite{M1}),
%%because $\sum_{k \geq 0} p_{(k)}(u) p_{(n-k)}(-u) = p_{(n)}(0) = \delta_{n,0}$:
It can also be deduced from the following relation 
(essentially the same as (2.9) in~\cite{M1}),
because $\sum_{k \geq 0} p_{(k)}(u) p_{(n-k)}(-u) = p_{(n)}(0) = \delta_{n,0}$:

\begin{theorem}\sl
%%   Assume the relations $c_0 = c'_0= 1$ and 
%%   $\sum_{k \geq 0} (-)^k c_k c'_{n-k} = \delta_{n,0}$
%%   between two sequences $\{ c_k \}_{k \geq 0}$ and $\{ c'_k \}_{k \geq 0}$.
   Assume that sequences $\{ c_k \}_{k \geq 0}$ and $\{ c'_k \}_{k \geq 0}$
   satisfy the relations $c_0 = c'_0= 1$ and 
   $\sum_{k \geq 0} (-)^k c_k c'_{n-k} = \delta_{n,0}$.
%%   Then we have
   Then, for two partitions $\lambda$ and $\mu$, we have
   $$
      \det 
      \left( 
	     c_{\lambda_i - \mu_j -i +j}
	  \right)_{1 \leq i,j \leq \operatorname{depth} \lambda}
	  = \det 
      \left( 
	     c'_{\lambda'_i - \mu'_j -i +j}
	  \right)_{1 \leq i,j \leq \operatorname{depth} \lambda'}.      
   $$
   Here we interpret $c_n$ and $c'_n$ as $0$ for $n<0$.
\end{theorem}
%
%%%%%%%%%%%%%%%%%%%%%%%%%%%%%%%%%%%%%%%%%%%%%%%%%%%%%%%%%%%%%%%%%%%%%%%%%%
%
\section{Expansions of $e$ and $h$}\label{sec:expansions_of_e_and_h}
%
%%%%%%%%%%%%%%%%%%%%%%%%%%%%%%%%%%%%%%%%%%%%%%%%%%%%%%%%%%%%%%%%%%%%%%%%%%
%

The expansion formulas for the functions $e_k$, $e^*_k$, $h_k$, and $h^*_k$ 
%%have their own aspects. 
have more interesting aspects. 
Some of these formulas are deduced 
from the results in the previous section as special cases,
but the others are not,
and we observe a mysterious duality in these formulas.
For simplicity, we introduce the following notation:
\begin{align*}
   h_{k}(x_1,\ldots,x_N ; u) &=  h_{k}(x_1 + u,\ldots,x_N + u), \\
   h^*_{k}(x_1,\ldots,x_N ; u) &=  h^*_{k}(x_1 + u,\ldots,x_N + u), \\
   e_{k}(x_1,\ldots,x_N ; u) &=  e_{k}(x_1 - u,\ldots,x_N - u), \\
   e^*_{k}(x_1,\ldots,x_N ; u) &=  e^*_{k}(x_1 - u,\ldots,x_N - u).
\end{align*}
These can be expanded as follows:
\begin{theorem}\label{thm:expansion_of_e}
   \sl For $k \geq 0$, we have
   \begin{align*}
      e_{k}(x_1,\ldots,x_N ; u)
	  & = \sum_{l \geq 0} {-N+k-1 \choose k-l} e_{l}(x_1,\ldots,x_N) p_{k-l}(u) \\
	  & = \sum_{l \geq 0} {-N+k-1 \choose k-l} e^*_{l}(x_1,\ldots,x_N) p^*_{k-l}(u), \\
	  e^*_{k}(x_1,\ldots,x_N ; u)
	  & = \sum_{l \geq 0} {-N+k-1 \choose k-l} e^*_{l}(x_1,\ldots,x_N) p_{k-l}(u).
   \end{align*}
\end{theorem}
\begin{theorem}\label{thm:expansion_of_h_1}
   \sl For $k \geq 0$, we have
   \begin{align*}
      h_{k}(x_1,\ldots,x_N ; u)
	  & = \sum_{l \geq 0} {N+k-1 \choose k-l} h_{l}(x_1,\ldots,x_N) p_{k-l}(u), \\
      h^*_{k}(x_1,\ldots,x_N ; u)
	  & = \sum_{l \geq 0} {N+k-1 \choose k-l} h^*_{l}(x_1,\ldots,x_N) p_{k-l}(u) \\
	  & = \sum_{l \geq 0} {N+k-1 \choose k-l} h_{l}(x_1,\ldots,x_N) p^*_{k-l}(u). 
   \end{align*}
\end{theorem}
Comparing these two theorems,
we observe a duality corresponding 
to the exchanges $e \leftrightarrow h^*$ and  $e^* \leftrightarrow h$. 
Note that the following are not equalities in general:
\begin{align*}
   e^*_{k}(x_1,\ldots,x_N ; u)
   & \ne \sum_{l \geq 0} {-N+k-1 \choose k-l} e_{l}(x_1,\ldots,x_N) p^*_{k-l}(u), \\
   h_{k}(x_1,\ldots,x_N ; u)
   & \ne \sum_{l \geq 0} {N+k-1 \choose k-l} h^*_{l}(x_1,\ldots,x_N) p^*_{k-l}(u).
\end{align*}
Theorem~\ref{thm:expansion_of_h_1} can be extended for negative integers $k$ as follows:
\begin{theorem}\label{thm:expansion_of_h_2}\sl
   For $k \in \mathbb{Z}$, we have
   \begin{align*}
      h_{k}(x_1,\ldots,x_N ; u)
	  & = \sum_{l \geq 0} {N+k-1 \choose N+l-1} h_{l}(x_1,\ldots,x_N) p_{k-l}(u), \\
      h^*_{k}(x_1,\ldots,x_N ; u)
	  & = \sum_{l \geq 0} {N+k-1 \choose N+l-1} h^*_{l}(x_1,\ldots,x_N) p_{k-l}(u) \\
	  & = \sum_{l \geq 0} {N+k-1 \choose N+l-1} h_{l}(x_1,\ldots,x_N) p^*_{k-l}(u). 
   \end{align*}
\end{theorem}
\begin{theorem}\label{thm:expansion_of_h_3}\sl
   For $k \in \mathbb{Z}$, we have
   \begin{align*}
      h_{k}(x_1,\ldots,x_N ; u)
	  & = \sum_{l \geq 0} {N+k-1 \choose l} h_{k-l}(x_1,\ldots,x_N) p_l(u), \\
      h^*_{k}(x_1,\ldots,x_N ; u)
	  & = \sum_{l \geq 0} {N+k-1 \choose l} h^*_{k-l}(x_1,\ldots,x_N) p_{l}(u) \\
	  & = \sum_{l \geq 0} {N+k-1 \choose l} h_{k-l}(x_1,\ldots,x_N) p^*_{l}(u). 
   \end{align*}
\end{theorem}
%
%%We can prove Theorems~\ref{thm:expansion_of_h_1}, 
%%\ref{thm:expansion_of_h_2}, and \ref{thm:expansion_of_h_3} 
We can prove Theorems~\ref{thm:expansion_of_h_1}--\ref{thm:expansion_of_h_3} 
easily in a way similar to the general expansion Theorem~\ref{thm:expansion_of_s}.
The proof of Theorem~\ref{thm:expansion_of_e} is as follows:

\medskip

\noindent
{\it Proof of Theorem~{\sl \ref{thm:expansion_of_e}}.}
First we prove the case of $k=N$.
Noting (\ref{eq:difference_product}), we have
\begin{align*}
   & \tilde{s}_{(1^N)}(x_1-u,\ldots,x_N-u) \\
   & \qquad = \det
   \begin{pmatrix}
      p_N(x_1-u) & \hdots & p_N(x_N-u) \\
      \vdots & & \vdots \\
      p_1(x_1-u) & \hdots & p_1(x_N-u) \\
   \end{pmatrix}\displaybreak[0]\\
   & \qquad = 
   (x_1-u) \cdots (x_N-u)
   \det
   \begin{pmatrix}
      p^*_{N-1}(x_1-u) & \hdots & p^*_{N-1}(x_N-u) \\
      \vdots & & \vdots \\
      p^*_0(x_1-u) & \hdots & p^*_0(x_N-u) \\
   \end{pmatrix}\displaybreak[0]\\
   & \qquad = (x_1-u) \cdots (x_N-u) \Delta(x_1-u,\ldots,x_N-u) \displaybreak[0]\\
   & \qquad = (x_1-u) \cdots (x_N-u) \Delta(x_1,\ldots,x_N) \displaybreak[0]\\
   & \qquad = \Delta(x_1,\ldots,x_N,u) \displaybreak[0]\\
   & \qquad 
   = \det
   \begin{pmatrix}
      p_N(x_1) & \hdots & p_N(x_N) & p_N(u)\\
      \vdots & & \vdots & \vdots \\
      p_0(x_1) & \hdots & p_0(x_N) & p_0(u) \\
   \end{pmatrix}.
\end{align*}
%%Considering the cofactor expansion along the last column,
By the cofactor expansion along the last column,
we see that this is equal to
$$
   \sum_{k=0}^N (-)^k \tilde{s}_{(1^{N-k})}(x_1,\ldots,x_N) p_k(u).
$$
Dividing both sides by the difference product $\Delta(x_1 - u,\ldots,x_N - u) = \Delta(x_1,\ldots,x_N)$,
we have
$$
   e_N(x_1-u,\ldots,x_N-u)
   = \sum_{k \geq 0} (-)^k p_{k}(u) e_{N-k}(x_1,\ldots,x_N). 
$$
Thus we have the assertion for $k = N$.

The other cases are deduced from this.
Indeed, on one hand, we have
$$
   e_N(x_1-u-w,\ldots,x_N-u-w)
   = \sum_{l \geq 0} (-)^l p_{l}(w) e_{N-l}(x_1-u,\ldots,x_N-u), 
$$
and on the other hand
\begin{align*}
   e_N(x_1-u-w,\ldots,x_N-u-w)
   & = \sum_{k \geq 0} (-)^k p_{k}(u+w) e_{N-k}(x_1,\ldots,x_N) \\
   & = \sum_{k \geq 0} \sum_{l \geq 0} (-)^k {k \choose l} p_{l}(w) p_{k-l}(u) e_{N-k}(x_1,\ldots,x_N).
\end{align*}
Comparing the coefficients of $p_{l}(w)$, 
we have the general case. 
\hfil\qed

\medskip

%%The following relation for the delta operator is easy from these expansions:
The following relation for the delta operator is easy to deduce from these expansions:
\begin{corollary}\sl
   We have
   \begin{align*}
     Q_u h_{k}(x_1,\ldots,x_N ; u)
	  & = (N+k-1) h_{k-1}(x_1,\ldots,x_N ; u), \\
	 Q_u h^*_{k}(x_1,\ldots,x_N ; u)
	  & = (N+k-1) h^*_{k-1}(x_1,\ldots,x_N ; u), \\
     Q_u e_{k}(x_1,\ldots,x_N ; u)
	  & = (-N+k-1) e_{k-1}(x_1,\ldots,x_N ; u), \\
     Q_u e^*_{k}(x_1,\ldots,x_N ; u)
	  & = (-N+k-1) e^*_{k-1}(x_1,\ldots,x_N ; u).
   \end{align*}
\end{corollary}
%
%%%%%%%%%%%%%%%%%%%%%%%%%%%%%%%%%%%%%%%%%%%%%%%%%%%%%%%%%%%%%%%%%%%%%%%%%%
%
\section{Generating functions}\label{sec:generating_functions}
%
%%%%%%%%%%%%%%%%%%%%%%%%%%%%%%%%%%%%%%%%%%%%%%%%%%%%%%%%%%%%%%%%%%%%%%%%%%
%
Combining Proposition~\ref{prop:special_e_and_h*} with the relations in the previous section,
we have the following relations
(put $k=N$ in Theorem~\ref{thm:expansion_of_e}
and $k=-N$ in Theorems~\ref{thm:expansion_of_h_2} and~\ref{thm:expansion_of_h_3}):
\begin{theorem}\label{thm:gen_fn_of_e}\sl
   We have
   \begin{align*}
      (u - x_1) \cdots (u - x_N)
      & = \sum_{l \geq 0} (-)^{l} e_l(x_1,\ldots,x_N) p_{N-l}(u) \\
      & = \sum_{l \geq 0} (-)^{l} e^*_l(x_1,\ldots,x_N) p^*_{N-l}(u).
   \end{align*}
\end{theorem}
\begin{theorem}\label{thm:gen_fn_of_h_1}\sl
   We have
   \begin{align*}
      \frac{1}{(u + x_1) \cdots (u + x_N)}
      & = \sum_{l \geq 0} (-)^l h_l(x_1,\ldots,x_N) p^*_{-N-l}(u) \\
      & = \sum_{l \geq 0} (-)^l h^*_l(x_1,\ldots,x_N) p_{-N-l}(u).
   \end{align*}
\end{theorem}
\begin{theorem}\label{thm:gen_fn_of_h_2}\sl
   We have
   \begin{align*}
      \frac{1}{(u + x_1) \cdots (u + x_N)}
      & =  (-)^{N-1} \sum_{l \geq 0} (-)^{l} h_{-N-l}(x_1,\ldots,x_N) p^*_{l}(u) \\
      & =  (-)^{N-1} \sum_{l \geq 0} (-)^{l} h^*_{-N-l}(x_1,\ldots,x_N) p_{l}(u).
   \end{align*}
\end{theorem}
We can regard the left hand sides of these equalities 
as ``generating functions'' of $e_k$, $e^*_k$, $h_k$, and $h^*_k$
represented as sums of multiples of $p_n(x)$ or $p^*_n(x)$
instead of the ordinary power.
\begin{remark}
%%Theorem~\ref{thm:gen_fn_of_e} also holds, 
The conclusion of Theorem~\ref{thm:gen_fn_of_e} also holds, 
even if $\{ p_n(u) \}$ is not of binomial type.
Namely, if $p_n(u)$ is a monic polynomial of degree $n$,
we have
$$
   (u - x_1) \cdots (u - x_N)
   = \sum_{l \geq 0} (-)^{l} e_l(x_1,\ldots,x_N) p_{N-l}(u).
$$
This is easily seen from the proof of Theorem~\ref{thm:expansion_of_e}.
%%Similarly, Theorems~\ref{thm:gen_fn_of_h_1} and \ref{thm:gen_fn_of_h_2} 
Similarly, the assertions of Theorems~\ref{thm:gen_fn_of_h_1} and \ref{thm:gen_fn_of_h_2} 
also hold, if $p_n(x)$ and $p^*_n(x)$ are monic polynomials of degree $n$
satisfying the relation
$$
   \frac{1}{x+y} 
   = \sum_{k \geq 0} (-)^k p_k(x) p^*_{-1-k}(y)
   = \sum_{k \geq 0} (-)^k p^*_k(x) p_{-1-k}(y).
$$
\end{remark}
%
%%%%%%%%%%%%%%%%%%%%%%%%%%%%%%%%%%%%%%%%%%%%%%%%%%%%%%%%%%%%%%%%%%%%%%%%%%
%
\section{Cauchy type relations}
%
%%%%%%%%%%%%%%%%%%%%%%%%%%%%%%%%%%%%%%%%%%%%%%%%%%%%%%%%%%%%%%%%%%%%%%%%%%
%
The relations in the previous section can be generalized as analogues
of the (dual) Cauchy identity.
In this section, 
%%we often abbreviate a function $f(x_1,\ldots,x_N)$ as $f(x)$ simply.
we often abbreviate a function $f(x_1,\ldots,x_N)$ simply as~$f(x)$.
\begin{theorem}\label{thm:Cauchy_1}\sl
   We have
%%   $$
%%      \prod_{1 \leq i \leq M} \prod_{1 \leq j \leq N} (y_j - x_i)
%%	  = \sum_{\lambda \supset \emptyset} 
%%	  (-)^{|\lambda|} s_{\lambda}(x) s_{\lambda^{\dagger}}(y) 
%%      = \sum_{\lambda \supset \emptyset} 
%%	  (-)^{|\lambda|} s^*_{\lambda}(x) s^*_{\lambda^{\dagger}}(y).
%%   $$
   $$
      \prod_{1 \leq i \leq M} \prod_{1 \leq j \leq N} (y_j - x_i)
	  = \sum_{\lambda}
	  (-)^{|\lambda|} s_{\lambda}(x) s_{\lambda^{\dagger}}(y) 
      = \sum_{\lambda}
	  (-)^{|\lambda|} s^*_{\lambda}(x) s^*_{\lambda^{\dagger}}(y).
   $$
   Here $\lambda$ runs over the Young diagrams satisfying 
   $\operatorname{depth}(\lambda) \leq M$ and
   $\operatorname{depth}(\lambda') \leq N$.
   Moreover we define $\lambda^{\dagger}$ by
   $$
      \lambda^{\dagger} = (N-\lambda_M,N-\lambda_{M-1},\ldots,N-\lambda_1)'.
   $$
\end{theorem}
\begin{theorem}\label{thm:Cauchy_2}\sl
   We have
%%   $$
%%      \prod_{1 \leq i \leq M} \prod_{1 \leq j \leq N} \frac{1}{y_j + x_i}
%%	  = \sum_{\lambda \supset \emptyset} (-)^{|\lambda|} 
%%	      s^*_{\lambda}(x) s_{\lambda^{\ddagger}}(y) 
%%      = \sum_{\lambda \supset \emptyset} (-)^{|\lambda|} 
%%	      s_{\lambda}(x) s^*_{\lambda^{\ddagger}}(y)
%%   $$
   $$
      \prod_{1 \leq i \leq M} \prod_{1 \leq j \leq N} \frac{1}{y_j + x_i}
	  = \sum_{\lambda} (-)^{|\lambda|} 
	      s^*_{\lambda}(x) s_{\lambda^{\ddagger}}(y) 
      = \sum_{\lambda} (-)^{|\lambda|} 
	      s_{\lambda}(x) s^*_{\lambda^{\ddagger}}(y)
   $$
   in $\mathbb{C}[x_1,\ldots,x_M]((y_1^{-1},\ldots,y_N^{-1}))$.
   Here $\lambda$ runs over the Young diagrams $\lambda$ satisfying
   $\operatorname{depth}(\lambda) \leq \operatorname{min}(M,N)$.
   Moreover we define $\lambda^{\ddagger}$ by
   $$
      \lambda^{\ddagger} 
	  = (-M-\lambda_{N},-M-\lambda_{N-1},\ldots,-M-\lambda_1).
   $$
\end{theorem}
From Theorem~\ref{thm:Cauchy_1} we see the following duality:
$$
   (-)^{|\lambda|} d_{\lambda\mu}(u) = (-)^{|\mu|} d_{\mu^{\dagger}\lambda^{\dagger}}(-u).
$$
Indeed, this follows by expanding 
$\prod_{1 \leq i \leq M} \prod_{1 \leq j \leq N} \frac{1}{y_j + x_i + u}$
in two ways.
Theorem~\ref{thm:duality_coeff} is immediate from this.
\begin{remark}
%%As in the previous section, 
%%Theorem~\ref{thm:Cauchy_1} holds
%%even if $\{ p_n(u) \}$ is not of binomial type.
%%Namely this holds if $p_n(u)$ is a monic polynomial of degree $n$.
%%Similarly Theorem~\ref{thm:Cauchy_2} also holds,
%%if $p_n(x)$ and $p^*_n(x)$ are monic polynomials of degree $n$ satisfying
As in the previous section, 
the conclusion of Theorem~\ref{thm:Cauchy_1} holds
even if $\{ p_n(u) \}$ is not of binomial type.
Namely it holds if $p_n(u)$ is a monic polynomial of degree $n$.
Similarly the conclusion of Theorem~\ref{thm:Cauchy_2} also holds
if $p_n(x)$ and $p^*_n(x)$ are monic polynomials of degree $n$ satisfying
$$
   \frac{1}{x+y} 
   = \sum_{k \geq 0} (-)^k p_k(x) p^*_{-1-k}(y)
   = \sum_{k \geq 0} (-)^k p^*_k(x) p_{-1-k}(y).
$$
\end{remark}
\begin{remark}
We can regard Theorem~\ref{thm:Cauchy_2} as 
a generalization of the following well-known relation
(the Cauchy identity):
%%$$
%%   \prod_{1 \leq i \leq M,\,1 \leq i \leq N}\frac{1}{1 - x_i y_j}
%%   = \sum_{\lambda \supset \emptyset} s^D_{\lambda}(x) s^D_{\lambda}(y).
%%$$
$$
   \prod_{1 \leq i \leq M,\,1 \leq i \leq N}\frac{1}{1 - x_i y_j}
   = \sum_{\lambda} s^D_{\lambda}(x) s^D_{\lambda}(y).
$$
\end{remark}

\noindent
{\it Proof of Theorem~{\sl \ref{thm:Cauchy_1}}.}
By (\ref{eq:difference_product}) we have
\begin{align*}
   & \Delta(y) \Delta(x) \cdot \prod_{1 \leq i \leq N} \prod_{1 \leq j \leq M} (y_j - x_i) 
   = \Delta(y_1,\ldots,y_N,x_1,\ldots,x_M) \\
   & \qquad 
   = \det
   \begin{pmatrix}
      p_{M+N-1}(y_1) & \hdots & p_{M+N-1}(y_N) & p_{M+N-1}(x_1) & \hdots & p_{M+N-1}(x_M) \\
      p_{M+N-2}(y_1) & \hdots & p_{M+N-2}(y_N) & p_{M+N-2}(x_1) & \hdots & p_{M+N-2}(x_M) \\
	  \vdots        &        & \vdots        & \vdots        &        & \vdots        \\
	  p_{0}(y_1)    & \hdots & p_{0}(y_N)    & p_{0}(x_1)    & \hdots & p_{0}(x_M)   
   \end{pmatrix}.
\end{align*}
Applying the Laplace expansion to the first $N$ columns, 
we see that this is equal to
%%$\sum_{\lambda \supset \emptyset}
%%(-)^{|\lambda|} \tilde{s}_{\lambda}(x) \tilde{s}_{\lambda^{\dagger}}(y)$.
$\sum_{\lambda}
(-)^{|\lambda|} \tilde{s}_{\lambda}(x) \tilde{s}_{\lambda^{\dagger}}(y)$.
Indeed we have 
$$
   \{ \lambda_i + M - i \,|\, 1 \leq i \leq M \} \cup 
   \{ \lambda^{\dagger}_j + N - j \,|\, 1 \leq j \leq N \} = \{ 0,1,\ldots, M+N-1 \}
$$
(recall (1.7) in \cite{M1}).
This means the first equality.
The second equality is similarly shown by replacing $p_k$ by $p^*_k$.
\hfill\qed

\medskip

To prove Theorem~\ref{thm:Cauchy_2}, we use the following well-known relation
%%(the left hand side is known as the ``Cauchy determinant''):
(the left hand side is known as the {\it Cauchy determinant}):
\begin{lemma}\label{lemma:Cauchy_det}\sl
   When $M=N$, we have
   $$
      \det \left( \frac{1}{x_i + y_j} \right)
	  = \frac{\Delta(x) \Delta(y)}{\prod_{1 \leq i,j \leq N} (x_i + y_j)}.
   $$
\end{lemma}
This relation is generalized as follows
(this can be proved by induction):
\begin{lemma}\label{lemma:gen_Cauchy_det}\sl
   When $N \geq M$,  we have
   $$
      \det\begin{pmatrix}
	  p^*_{N-M-1}(y_1)    & p^*_{N-M-1}(y_2)    & \hdots & p^*_{N-M-1}(y_N)    \\
	  \vdots              & \vdots              &        & \vdots              \\
	  p^*_1(y_1)          & p^*_1(y_2)          & \hdots & p^*_1(y_N)          \\
	  p^*_0(y_1)          & p^*_0(y_2)          & \hdots & p^*_0(y_N)          \\
      \frac{1}{x_1 + y_1} & \frac{1}{x_1 + y_2} & \hdots & \frac{1}{x_1 + y_N} \\
	  \frac{1}{x_2 + y_1} & \frac{1}{x_2 + y_2} & \hdots & \frac{1}{x_2 + y_N} \\
	  \vdots              & \vdots              &        & \vdots              \\
	  \frac{1}{x_M + y_1} & \frac{1}{x_M + y_2} & \hdots & \frac{1}{x_M + y_N} \\
      \end{pmatrix}
	  = \frac{\Delta(x) \Delta(y)}{\prod_{1 \leq i \leq M, \, 1 \leq j \leq N} (x_i + y_j)}.
   $$
\end{lemma}

\noindent
{\it Proof of Theorem~{\sl \ref{thm:Cauchy_2}}} ({\it the case of $M=N$}).
Using the Cauchy--Binet formula (Proposition~\ref{prop:Cauchy--Binet}), we have
\begin{align*}
   & \det \left( \frac{1}{x_i + y_j} \right)_{1 \leq i,j \leq N} \displaybreak[0]\\
   & \qquad = \det 
   \left( 
      \sum_{k \geq 0} (-)^k  p_k(x_i) p^*_{-1-k}(y_j)
   \right)_{1 \leq i,j \leq N} \displaybreak[0]\\
   & \qquad = \sum_{0 \leq k_1 < \cdots < k_N} 
   \det((-)^{k_i} p_{k_i}(x_j))_{1 \leq i,j \leq N} 
   \det(p^*_{-1-k_i}(y_j))_{1 \leq i,j \leq N} \displaybreak[0]\\
   & \qquad = \sum_{0 \leq k_1 < \cdots < k_N} (-)^{k_1 + \cdots + k_N}
   \det(p_{k_i}(x_j))_{1 \leq i,j \leq N} 
   \det(p^*_{-1-k_i}(y_j))_{1 \leq i,j \leq N} \displaybreak[0].
\end{align*}
%%Here, the first determinant in the right hand side is equal to
Here, the first determinant on the right hand side is equal to
$$
   \det(p_{\lambda_i+N-i}(x_j))_{1 \leq i,j \leq N} 
   = \tilde{s}_{\lambda}(x),
$$
where we define $\lambda_i$ by 
$$
   k_N = \lambda_1 + N -1, \quad
   k_{N-1} = \lambda_2 + N -2, \quad
   \ldots, \quad
   k_1 = \lambda_N + 0.
$$ 
On the other hand, the second determinant is equal to
$$
   \det(p^*_{-1-N-\lambda_i+i}(y_j))_{1 \leq i,j \leq N}
   = \tilde{s}^*_{\lambda^{\ddagger}}(y).
$$
Thus we have
%%$$
%%   \det \left( \frac{1}{x_i + y_j} \right)_{1 \leq i,j \leq N} 
%%   = \sum_{\lambda \supset \emptyset} (-)^{|\lambda|}
%%   \tilde{s}_{\lambda}(x) \tilde{s}^*_{\lambda^{\ddagger}}(y),
%%$$
$$
   \det \left( \frac{1}{x_i + y_j} \right)_{1 \leq i,j \leq N} 
   = \sum_{\lambda} (-)^{|\lambda|}
   \tilde{s}_{\lambda}(x) \tilde{s}^*_{\lambda^{\ddagger}}(y),
$$
and the assertion is immediate by dividing this by $\Delta(x)\Delta(y)$. 
\hfill\qed

\medskip

\noindent
{\it Proof of Theorem~{\sl \ref{thm:Cauchy_2}}} ({\it the case of $N > M$}).
%%Let us denote by $A$ the matrix in the left hand side of Lemma~\ref{lemma:gen_Cauchy_det}.
Let us denote by $A$ the matrix on the left hand side of Lemma~\ref{lemma:gen_Cauchy_det}.
This can be expressed as
$$
   A = 
   \begin{pmatrix}
      p^*_{N-M-1}(y_1) & \hdots & p^*_{N-M-1}(y_N) \\
	  \vdots              &        & \vdots          \\
	  p^*_{0}(y_1) & \hdots & p^*_{0}(y_N) \\
      \sum_k (-)^k p_k(x_1) p^*_{-1-k}(y_1) & \hdots & \sum_k (-)^k p_k(x_1) p^*_{-1-k}(y_N) \\
	  \vdots              &        & \vdots          \\
      \sum_k (-)^k p_k(x_M) p^*_{-1-k}(y_1) & \hdots & \sum_k (-)^k p_k(x_M) p^*_{-1-k}(y_N)
   \end{pmatrix}
   = BC
$$
with the $N \times \infty$ matrix $B$
and the $\infty \times N$ matrix $C$ defined by
\begin{align*}
   B & = 
   \begin{pmatrix}
   	  1 & \hdots & 0 & 0 & 0 & \hdots & \\
	  \vdots & \ddots & \vdots & \vdots & \vdots & &   \\
	  0 & \hdots & 1 & 0 & 0 & \hdots & \\
      0 & \hdots & 0 & (-)^0 p_0(x_1) & (-)^1 p_1(x_1) & \hdots & \\
	  \vdots &        & \vdots & \vdots & \vdots & & \\
	  0 & \hdots & 0 & (-)^0 p_0(x_M) & (-)^1 p_1(x_M) & \hdots & \\
   \end{pmatrix}, \\ 
   C & = 
   \begin{pmatrix}
      p^*_{N-M-1}(y_1) & \hdots & p^*_{N-M-1}(y_N) \\
	  p^*_{N-M-2}(y_1) & \hdots & p^*_{N-M-2}(y_N) \\
	  \vdots & & \vdots    
   \end{pmatrix}. 
\end{align*}
Applying the Cauchy--Binet formula (Proposition~\ref{prop:Cauchy--Binet}) 
to this relation $A=BC$, we have
$$
   \det A
   = \det BC
   = \sum_{1 \leq i_1 < \cdots < i_N} \det B_{(1,\ldots,N), (i_1,\ldots,i_N)} \det C_{(i_1,\ldots,i_N), (1,\ldots,N)}.
$$
Note that $\det B_{(1,\ldots,N), (i_1,\ldots,i_N)} = 0$,
unless $(i_1,\ldots,i_{N-M}) = (1,\ldots,N-M)$.
%%Thus this is equal to
Thus we have
\begin{align*}
   \det A =
   & \sum_{0 \leq k_1 < \cdots < k_M}
   \det
   \begin{pmatrix}
   	  1 & \hdots & 0 & 0 & \hdots &  0 \\
	  \vdots & \ddots & \vdots & \vdots & & \vdots \\
	  0 & \hdots & 1 & 0 & \hdots & 0 \\
      0 & \hdots & 0 & (-)^{k_1} p_{k_1}(x_1) & \hdots & (-)^{k_M} p_{k_M}(x_1) \\
	  \vdots &        & \vdots & \vdots & & \vdots \\
	  0 & \hdots & 0 & (-)^{k_1} p_{k_1}(x_M) &\hdots & (-)^{k_M} p_{k_M}(x_M) \\
   \end{pmatrix} \\
   & \qquad
   \cdot
   \det \begin{pmatrix}
      p^*_{N-M-1}(y_1) & \hdots & p^*_{N-M-1}(y_N) \\
 	  \vdots & & \vdots \\
      p^*_{0}(y_1) & \hdots & p^*_{0}(y_N) \\
	  p^*_{-1-k_1}(y_1) & \hdots & p^*_{-1-k_1}(y_N) \\
 	  \vdots & & \vdots \\
      p^*_{-1-k_M}(y_1) & \hdots & p^*_{-1-k_M}(y_N)
   \end{pmatrix}.
\end{align*}
%%On one hand, the first determinant is equal to
The first determinant is equal to
\begin{align*}
   & \det
   \begin{pmatrix}
      (-)^{k_1} p_{k_1}(x_1) & \hdots & (-)^{k_M} p_{k_M}(x_1) \\
	  \vdots & & \vdots \\
	  (-)^{k_1} p_{k_1}(x_M) &\hdots  & (-)^{k_M} p_{k_M}(x_M) \\
   \end{pmatrix} \displaybreak[0]\\
   & \qquad 
   =  (-)^{k_1 + \cdots + k_M}
   \det
   \begin{pmatrix}
      p_{k_1}(x_1) & \hdots & p_{k_M}(x_1) \\
	  \vdots & & \vdots \\
	  p_{k_1}(x_M) & \hdots & p_{k_M}(x_M) \\
   \end{pmatrix} \displaybreak[0]\\
   & \qquad
   =  (-)^{k_1 + \cdots + k_M}
   (-)^{\frac{M(M-1)}{2}}
   \det
   \begin{pmatrix}
      p_{k_M}(x_1) & \hdots & p_{k_1}(x_1) \\
	  \vdots & & \vdots \\
	  p_{k_M}(x_M) & \hdots & p_{k_1}(x_M) \\
   \end{pmatrix} \displaybreak[0]\\
   & \qquad
   =  (-)^{\lambda_1 + \cdots + \lambda_M}
   \det
   \begin{pmatrix}
      p_{\lambda_1 + M-1}(x_1) & \hdots & p_{\lambda_M}(x_1) \\
	  \vdots & & \vdots \\
	  p_{\lambda_1 + M-1}(x_M) & \hdots & p_{\lambda_M}(x_M) \\
   \end{pmatrix} \displaybreak[0]\\
   & \qquad
   = (-)^{|\lambda|} \tilde{s}_{\lambda}(x).
\end{align*}
Here we define $\lambda_1,\ldots,\lambda_M$
by 
$$
   k_M = \lambda_1 + M-1,\quad
   k_{M-1} = \lambda_2 + M-2,\quad
   \ldots,\quad
   k_1 = \lambda_M + 0.
$$
%%On the other hand, the second determinant is equal to
The second determinant in the formula for $\det A$ is equal to
$$
   \det \begin{pmatrix}
      p^*_{\lambda^{\ddagger}_1 + N-1}(y_1) & \hdots & p^*_{\lambda^{\ddagger}_1 + N-1}(y_N) \\
	  p^*_{\lambda^{\ddagger}_2 + N-2}(y_1) & \hdots & p^*_{\lambda^{\ddagger}_2 + N-2}(y_N) \\
 	  \vdots & & \vdots \\
      p^*_{\lambda^{\ddagger}_{N} + 0}(y_1) & \hdots & p^*_{\lambda^{\ddagger}_{N} + 0}(y_N) \\
   \end{pmatrix}
   = \tilde{s}^*_{\lambda^{\ddagger}}(y).
$$
Here we put
%%$$
%%   \lambda^{\ddagger}_1 = -M-\lambda_M,\,
%%   \lambda^{\ddagger}_2 = -M-\lambda_{M-1},\,
%%   \ldots,\,
%%   \lambda^{\ddagger}_N = -M-\lambda_1.
%%$$
$$
   \lambda^{\ddagger}_1 = -M-\lambda_M,\quad
   \lambda^{\ddagger}_2 = -M-\lambda_{M-1},\quad
   \ldots,\quad
   \lambda^{\ddagger}_N = -M-\lambda_1.
$$
Thus we have
%%$$
%%   \det A
%%   = \sum_{\lambda \supset \emptyset} (-)^{|\lambda|}
%%   \tilde{s}_{\lambda}(x) \tilde{s}^*_{\lambda^{\ddagger}}(y).
%%$$
$$
   \det A
   = \sum_{\lambda} (-)^{|\lambda|}
   \tilde{s}_{\lambda}(x) \tilde{s}^*_{\lambda^{\ddagger}}(y).
$$
%%This means our assertion.
This yields our assertion.
\hfil\qed

\medskip

%%The proof in the case $N < M$ is almost the same, so that we omit it.
The proof in the case $N < M$ is almost the same, so we omit it.

%
%%%%%%%%%%%%%%%%%%%%%%%%%%%%%%%%%%%%%%%%%%%%%%%%%%%%%%%%%%%%%%%%%%%%%%%%%%
%
\section{Capelli type elements}\label{sec:capelli_type_elements}
%
%%%%%%%%%%%%%%%%%%%%%%%%%%%%%%%%%%%%%%%%%%%%%%%%%%%%%%%%%%%%%%%%%%%%%%%%%%
%
Our Schur type functions are useful to express
the eigenvalues of Capelli type central elements 
of the universal enveloping algebras of the classical Lie algebras.
Before stating this, 
we recall these Capelli type elements in this section.

These central elements have been investigated 
in the study of Capelli type identities.
See \cite{HU}, \cite{MN}, \cite{O},  [U1--5], \cite{IU}, and [I1--6] for 
the Capelli identity and its generalizations.

%%%%%%%%%%%%%%%%%%%%%%%%%%%%%%%%%%%%%%%%%%%%%%%%%%%%%%%%%%%%%%%%%%%%%
\subsection{}
First, we recall the Capelli elements,
famous central elements of the universal enveloping algebra $U(\mathfrak{gl}_N)$
of the general linear Lie algebra $\mathfrak{gl}_N$.
Let $E_{ij}$ be the standard basis of $\mathfrak{gl}_N$,
and consider the matrix 
$E = (E_{ij})_{1 \leq i,j \leq N}$ in 
$\operatorname{Mat}_N(U(\mathfrak{gl}_N))$.
Then the following determinant
%%is known as the ``Capelli determinant'' (\cite{Ca1}, \cite{HU}, \cite{U1}):
is known as the {\it Capelli element} (\cite{Ca1}, \cite{HU}, \cite{U1}):
$$
   C^{\mathfrak{gl}_N}(u) = \det(E - u \bold{1} + \operatorname{diag} \natural_N).
$$
%%Here $\bold{1}$ means the unit matrix,
%%and $\natural_N$ means the sequence $\natural_N = (N-1,N-2,\ldots,0)$ of length $N$.
Here $\bold{1}$ is the unit matrix,
and $\natural_N$ is the sequence $\natural_N = (N-1,N-2,\ldots,0)$ of length $N$.
Moreover ``$\det$'' means a non-commutative determinant
%%called the ``column-determinant.''
called the {\it column-determinant.}
Namely, for a square matrix $Z = (Z_{ij})$ whose entries are non-commutative,
we put
$$
   \det Z = \sum_{\sigma \in \mathfrak{S}_N} \operatorname{sgn}(\sigma) 
   Z_{\sigma(1)1} Z_{\sigma(2)2} \cdots Z_{\sigma(N)N}.
$$
%%This Capelli element $C^{\frak{gl}_N}(u)$ is known to be central in $U(\mathfrak{gl}_N)$.
The Capelli element $C^{\frak{gl}_N}(u)$ is known to be central in $U(\mathfrak{gl}_N)$.

This is generalized to the sums of minors
\begin{align}\label{eq:definition_of_C}
   C^{\mathfrak{gl}_N}_k(u)
   = \sum_{1 \leq i_1 < \cdots < i_k \leq N} 
     \det(E_I - u \bold{1} + \operatorname{diag} \natural_k).
\end{align}
Here we put $Z_I = (Z_{i_a i_b})_{1 \leq a,b \leq k}$
for $I = (i_1,\ldots,i_k)$ and $Z =(Z_{ij})$.
%%We call this the Capelli element of degree $k$.
We call this the {\it Capelli element of degree $k$.}

Moreover we can consider the following analogue using the permanent \cite{N}:
\begin{align}\label{eq:definition_of_D}
   D^{\mathfrak{gl}_N}_k(u)
   = \sum_{1 \leq i_1 \leq \cdots \leq i_k \leq N}
   \frac{1}{I!} 
   \operatorname{per}(E_I + u \bold{1}_I - \bold{1}_I \operatorname{diag} \natural_k).
\end{align}
%%Here ``$\operatorname{per}$'' means the column-permanent.
Here ``$\operatorname{per}$'' means the {\it column-permanent.}
%%Namely for a square matrix $Z = (Z_{ij})$ of size $N$, we put
Namely, for a square matrix $Z = (Z_{ij})$ of size $N$, we put
$$
   \operatorname{per} Z = \sum_{\sigma\in\mathfrak{S}_N} Z_{\sigma(1)1} \cdots Z_{\sigma(N)N}.
$$
Here we put $I! = m_1! \cdots m_N!$,
%%where $m_1, \ldots, m_N$ mean the multiplicities of $I = (i_1,\ldots,i_k)$:
where $m_1, \ldots, m_N$ are the multiplicities of $I = (i_1,\ldots,i_k)$:
$$
    I = (i_1,\ldots,i_k) 
      = (\overbrace{1,\ldots,1}^{m_1},
         \overbrace{2,\ldots,2}^{m_2},\ldots,
         \overbrace{N,\ldots,N}^{m_N}).
$$
Since $I$ has some multiplicities in general, 
$Z_I = (Z_{i_a i_b})_{1 \leq a,b \leq k}$ is not a submatrix of $Z$ necessarily.

%%These $C^{\mathfrak{gl}_N}_k(u)$ and $D^{\mathfrak{gl}_N}_k(u)$
The elements $C^{\mathfrak{gl}_N}_k(u)$ and $D^{\mathfrak{gl}_N}_k(u)$
are also central in the universal enveloping algebra
(actually these are generators of the center of the universal enveloping algebra;
see \cite{HU}, \cite{I4}, \cite{N}, \cite{U1} for the details).

\begin{theorem} \sl
   For any $u \in \mathbb{C}$,
   $C^{\mathfrak{gl}_N}_k(u)$ and $D^{\mathfrak{gl}_N}_k(u)$ are central in $U(\mathfrak{gl}_N)$.
\end{theorem}

These central elements act on the irreducible representations 
as scalar operators by Schur's lemma.
These values (the eigenvalues) can be calculated 
by noting the following triangular decomposition:
$$
   \mathfrak{gl}_N = \mathfrak{n}^- \oplus \mathfrak{h} \oplus \mathfrak{n}^+.
$$
Here $\mathfrak{n}^-$, $\mathfrak{h}$, and $\mathfrak{n}^+$
are the subalgebras of $\mathfrak{gl}_N$ spanned by the elements 
$F^{\mathfrak{gl}_N}_{ij}$ such that
$i>j$, $i=j$, and $i<j$, respectively.
Namely, the entries in the lower triangular part,
in the diagonal part, and in the upper triangular part  
of the matrix $E^{\mathfrak{gl}_N}$
belong to $\mathfrak{n}^-$, $\mathfrak{h}$, and $\mathfrak{n}^+$, respectively.

For example, we can calculate the eigenvalue of $C^{\mathfrak{gl}_N}(u)$ as follows.
Let $\pi$ be the irreducible representation
determined by the partition $(\lambda_1,\ldots,\lambda_N)$,
and consider the action of $C^{\mathfrak{gl}_N}(u)$ 
to the highest weight vector $v$:
$$
   \pi(C^{\mathfrak{gl}_N}(u)) v
   = \sum_{\sigma \in \mathfrak{S}_N} \operatorname{sgn}(\sigma)
   \pi(E_{\sigma(1)1}(-u+N-1)) \cdots \pi(E_{\sigma(N)N}(-u+0)) v.
$$
Then, 
among these $N!$ terms,
there remains
only one term corresponding to 
$\sigma = 
\left(\begin{smallmatrix} 
1 & 2 & \hdots & N \\ 
1 & 2 & \hdots & N 
\end{smallmatrix}\right)$,
and the other $N! - 1$ terms all vanish,
because $\pi(E_{ij}) v = 0$ for $i < j$
and $\pi(E_{ii}) v = \lambda_i v$.
Thus we have
\begin{align*}
   \pi(C^{\mathfrak{gl}_N}(u)) v
   & = \pi(E_{11}(-u+N-1)) \cdots \pi(E_{NN}(-u+0)) v \\
   & = (\lambda_1 -u+N-1) \cdots (\lambda_N -u+0) v.
\end{align*}
and we see that $C^{\mathfrak{gl}_N}(u)$
%%acts as the multiplication by 
acts as multiplication by 
the scalar $(\lambda_1 -u+N-1) \cdots (\lambda_N -u+0)$ on $\pi$.

%%We can write down the eigenvalue of $C^{\mathfrak{gl}_N}_k(u)$ similarly
%%by considering the action of each column-determinant 
%%in the right hand side of (\ref{eq:definition_of_C}) to the highest weight vector $v$.
%%Indeed, also in this case,
%%only one term remains
%%among $k!$ terms in each column-determinant.

We can write down the eigenvalue of $C^{\mathfrak{gl}_N}_k(u)$ similarly
by considering the action to the highest weight vector $v$
of each column-determinant on the right hand side of (\ref{eq:definition_of_C}).
Indeed, also in this case,
only one term remains of the $k!$ terms in each column-determinant.

The calculation of the eigenvalue of $D^{\mathfrak{gl}_N}_k(u)$
is essentially same (but a bit more complicated).
For this, we consider the action of each column-permanent 
%%in the right hand side of (\ref{eq:definition_of_D}) to $v$.
on the right hand side of (\ref{eq:definition_of_D}) to $v$.
%%Then, among the $k!$ terms in each column-permanent,
Then, of the $k!$ terms in each column-permanent,
there remain only $I!$ terms 
corresponding to $\sigma$ such that $i_{\sigma(a)} = i_a$ for $a = 1,\ldots,k$,
and these $I!$ terms are all equal to a scalar multiple of $v$.

In this way, we can write down
the eigenvalues of $C^{\mathfrak{gl}_N}_k(u)$ and $D^{\mathfrak{gl}_N}_k(u)$.
However the results are not so simple.
As seen in Section~\ref{sec:eigenvalues} below,
these eigenvalues are actually expressed by using the factorial (shifted) Schur functions.

\begin{remark}
%%We note some preceding results:
We note some previous results:

\smallskip

\noindent
(1)
We can also express
the elements $C^{\mathfrak{gl}_N}_k(u)$ and $D^{\mathfrak{gl}_N}_k(u)$ 
in terms of the ``symmetrized determinant'' and the ``symmetrized permanents''
(\cite{IU}, [I1--6]).
%%Under these expressions, 
From these expressions, 
we can easily see the centrality of these elements.

\smallskip

\noindent
(2)
In \cite{O},
Okounkov introduced a class of central elements of $U(\mathfrak{gl}_N)$,
%%and called them the ``quantum immanants.'' 
and called them the {\it quantum immanants.} 
These elements $\mathbb{S}_{\mu}$ indexed by partitions $\mu$
%%are expressed in terms of a determinant type function called ``immanant,''
are expressed in terms of a determinant type function called {\it immanant,}
and form 
a basis of the center of $U(\mathfrak{gl}_N)$ as a vector.
The central elements $C^{\mathfrak{gl}_N}_k(0)$ and $D^{\mathfrak{gl}_N}_k(0)$ 
can be regarded as the cases of $\mu = (1^k)$ and $\mu = (k)$:
\begin{align}\label{eq:relation_between_C_D_and_QI}
   C^{\mathfrak{gl}_N}_k(0) = \mathbb{S}_{(1^k)}, \qquad
   D^{\mathfrak{gl}_N}_k(0) = \mathbb{S}_{(k)}.
\end{align}
Okounkov also gave a generalization of the Capelli identity 
%%(called ``higher Capelli identities'')
(called {\it higher Capelli identities})
for the quantum immanants.
In Section~\ref{sec:eigenvalues},
we will see the eigenvalues of the quantum immanants
together with those of $C^{\mathfrak{gl}_N}_k(u)$ and $D^{\mathfrak{gl}_N}_k(u)$.
\end{remark}

%%%%%%%%%%%%%%%%%%%%%%%%%%%%%%%%%%%%%%%%%%%%%%%%%%%%%%%%%%%%%%%%%%%%%
\subsection{}
Next, we recall analogues of the Capelli elements
in the universal enveloping algebras of the orthogonal Lie algebras given in \cite{W}.

%%Let $S \in \operatorname{Mat}_N(\mathbb{C})$ be 
%%a non-degenerate symmetric matrix of size $N$. 
Let $S \in \operatorname{Mat}_N(\mathbb{C})$ be 
a non-singular symmetric matrix of size $N$. 
We can realize the orthogonal Lie group as the isometry group 
with respect to the bilinear form determined by $S$:
$$
   O(S) = \{ g \in GL_N \,|\, {}^t\!g S g = S \}.
$$
The corresponding Lie algebra is expressed as 
$$
   \mathfrak{o}(S) = \{ Z \in \mathfrak{gl}_N \,|\, {}^t\!Z S + S Z = 0\}.
$$
As generators of this $\mathfrak{o}(S)$,
we can take $F^{\mathfrak{o}(S)}_{ij} = E_{ij} - S^{-1} E_{ji} S$,
where $E_{ij}$ is the standard basis of $\frak{gl}_N$.
We introduce the $N \times N$ matrix $F^{\mathfrak{o}(S)}$
whose $(i,j)$th entry is this generator:
$F^{\mathfrak{o}(S)} = (F^{\mathfrak{o}(S)}_{ij})_{1 \leq i,j \leq N}$.
We regard this matrix as an element of $\operatorname{Mat}_N(U(\mathfrak{o}(S)))$.

In particular, in the case of $S = S_0 = (\delta_{i,N+1-j})$,
the corresponding orthogonal Lie algebra is expressed as follows:
$$
   \mathfrak{o}(S_0) 
   = \{ Z = (Z_{ij}) \in \mathfrak{gl}_N \,|\, Z_{ij} + Z_{N+1-j, N+1-i} = 0 \}.
$$
%%We call this the ``split realization'' of the orthogonal Lie algebra.
We call this the {\it split realization} of the orthogonal Lie algebra.
In this case, we can take the following triangular decomposition:
$$
   \mathfrak{o}(S_0) = \mathfrak{n}^- \oplus \mathfrak{h} \oplus \mathfrak{n}^+.
$$
Here $\mathfrak{n}^-$, $\mathfrak{h}$, and $\mathfrak{n}^+$
are the subalgebras of $\mathfrak{o}(S_0)$ spanned by the elements 
$F^{\mathfrak{o}(S_0)}_{ij}$ such that
$i>j$, $i=j$, and $i<j$, respectively.
Namely, the entries in the lower triangular part,
in the diagonal part, and in the upper triangular part  
of the matrix $E^{\mathfrak{gl}_N}$
belong to $\mathfrak{n}^-$, $\mathfrak{h}$, and $\mathfrak{n}^+$, respectively.
Thus, if there is a central element of $U(\frak{o}(S_0))$
%%which is expressed as ``the column-determinant of $F^{\mathfrak{o}(S_0)}$,''
which is expressed as the column-determinant of $F^{\mathfrak{o}(S_0)}$,
we can easily calculate its eigenvalue 
in a way similar to the case of $C^{\mathfrak{gl}_N}(u)$.
In fact,  the following central element was given by Wachi \cite{W}:

\begin{theorem}[Wachi] \sl
   For any $u \in \mathbb{C}$, 
%%   the following is central in $U(\mathfrak{o}(S_0))$:
   the following element is central in $U(\mathfrak{o}(S_0))$:
   $$
      C^{\mathfrak{o}_N}(u)
      = \det(F^{\mathfrak{o}(S_0)} - u \bold{1} 
        + \operatorname{diag} \tilde{\natural}_N).
   $$
   Here $\tilde{\natural}_N$ is the following sequence of length $N$:
   $$
      \tilde{\natural}_N =
      \begin{cases}
         (\tfrac{N}{2}-1, \tfrac{N}{2}-2, \ldots, 0, 
             0, \ldots, -\tfrac{N}{2}+1), & \text{$N$: even}, \\
         (\tfrac{N}{2}-1, \tfrac{N}{2}-2, \ldots, \tfrac{1}{2}, 0, 
             -\tfrac{1}{2}, \ldots, -\tfrac{N}{2}+1), & \text{$N$: odd}.
      \end{cases}
   $$
\end{theorem}

This can be generalized as follows:

\begin{theorem}[Wachi] \sl
   For any $u \in \mathbb{C}$, 
%%   the following is central in $U(\mathfrak{o}(S_0))$:
   the following element is central in $U(\mathfrak{o}(S_0))$:
   $$
      C^{\mathfrak{o}_N}_k(u)
      = \sum_{1 \leq i_1 < \cdots < i_k \leq N}
        \det(\widetilde{F}^{\frak{o}(S_0)}_I - u \bold{1} 
        + \operatorname{diag}(\tfrac{k}{2}-1,\tfrac{k}{2}-2,\ldots,-\tfrac{k}{2})).
   $$
   Here $\widetilde{F}^{\mathfrak{o}(S_0)}$ is defined by
   $$
      \widetilde{F}^{\mathfrak{o}(S_0)} = 
      \begin{cases}
         F^{\mathfrak{o}(S_0)} + \operatorname{diag}(0,\ldots,0,1,\ldots,1), 
         & \qquad \text{$N$: even}, \\
         F^{\mathfrak{o}(S_0)} + \operatorname{diag}(0,\ldots,0,\tfrac{1}{2},1,\ldots,1), 
         & \qquad \text{$N$: odd}.
      \end{cases}
   $$
   Here the numbers of $0$'s and $1$'s are equal to $[N/2]$.
\end{theorem}

In Section~\ref{sec:eigenvalues}, 
we will express the eigenvalues of these central elements
in terms of the Schur type functions
associated with the central difference.

%%%%%%%%%%%%%%%%%%%%%%%%%%%%%%%%%%%%%%%%%%%%%%%%%%%%%%%%%%%%%%%%%%%%%
\subsection{}
Similarly, we can construct central elements in the universal enveloping algebra
of the symplectic Lie algebra.
However, these are expressed in terms of permanents (not in terms of determinants).

%%Let $J \in \operatorname{Mat}_N(\mathbb{C})$ be 
%%a non-degenerate alternating matrix of size $N$. 
Let $J \in \operatorname{Mat}_N(\mathbb{C})$ be 
a non-singular alternating matrix of size $N$. 
We can realize the symplectic Lie group as the isometry group 
with respect to the bilinear form determined by $J$:
$$
   Sp(J) = \{ g \in GL_N \,|\, {}^t\!g J g = J \}.
$$
The corresponding Lie algebra is expressed as
$$
   \mathfrak{sp}(J) = \{ Z \in \mathfrak{gl}_N \,|\, {}^t\!Z J + J Z = 0\}.
$$
As generators of this $\mathfrak{sp}(J)$,
we can take $F^{\mathfrak{sp}(J)}_{ij} = E_{ij} - J^{-1} E_{ji} J$,
where $E_{ij}$ is the standard basis of $\frak{gl}_N$.
We introduce the $N \times N$ matrix $F^{\mathfrak{sp}(J)}$
whose $(i,j)$th entry is this generator:
$F^{\mathfrak{sp}(J)} = (F^{\mathfrak{sp}(J)}_{ij})_{1 \leq i,j \leq N}$.
We regard this matrix as an element of $\operatorname{Mat}_N(U(\mathfrak{sp}(J)))$.

\makeatletter
\def\iddots{\mathinner{\mkern1mu\raise1\p@
    \hbox{.}\mkern2mu\raise5\p@\hbox{.}\mkern2mu
    \raise9\p@\vbox{\kern7\p@\hbox{.}}\mkern1mu}}
\makeatother

%%We consider the split realization. 
We consider the split realization of the symplectic Lie algebra. 
Namely we consider the case of 
$$
   J = J_0 =
   \left(
   \begin{smallmatrix}
      & & & & & 1 \\
      & & & & \iddots & \\
      & & & 1 & & \\
      & & -1 & & & \\
      & \iddots & & & & \\
      -1 & & & & &
   \end{smallmatrix}
   \right).
$$
In this case, we can take the following triangular decomposition:
$$
   \mathfrak{sp}(J_0) = \mathfrak{n}^- \oplus \mathfrak{h} \oplus \mathfrak{n}^+.
$$
Here $\mathfrak{n}^-$, $\mathfrak{h}$, and $\mathfrak{n}^+$
are the subalgebras of $\mathfrak{sp}(J_0)$ spanned by the elements 
$F^{\mathfrak{sp}(J_0)}_{ij}$ such that
$i>j$, $i=j$, and $i<j$, respectively.
Namely, the entries in the lower triangular part,
in the diagonal part, and in the upper triangular part  
of the matrix $E^{\mathfrak{gl}_N}$
belong to $\mathfrak{n}^-$, $\mathfrak{h}$, and $\mathfrak{n}^+$, respectively.
In this case, we can construct central elements of the universal enveloping algebra
using the column-permanent \cite{I6}:

\begin{theorem}[Itoh] \sl
   For any $u \in \mathbb{C}$, 
%%   the following is central in $U(\mathfrak{sp}(J_0))$:
   the following element is central in $U(\mathfrak{sp}(J_0))$:
   $$
      D^{\mathfrak{sp}_N}_k(u) 
      = \sum_{1 \leq i_1 \leq \cdots \leq i_k \leq N} 
        \frac{1}{I!} \operatorname{per}(\widetilde{F}^{\mathfrak{sp}(J_0)}_I
      + u \bold{1}_I 
      - \bold{1}_I 
        \operatorname{diag}(\tfrac{k}{2}-1,\tfrac{k}{2}-2,\ldots,-\tfrac{k}{2})).
   $$
   Here we put
   $$
      \widetilde{F}^{\mathfrak{sp}(J_0)} 
      = F^{\mathfrak{sp}(J_0)} - \operatorname{diag}(0,\ldots,0,1,\ldots,1),
   $$
   where the numbers of $0$'s and $1$'s are equal to $N/2$.
\end{theorem}

We can easily calculate the eigenvalue of this central element
noting the triangular decomposition of $\mathfrak{sp}(J_0)$
and the definition of the column-permanent.
In Section~\ref{sec:eigenvalues}, 
we will express them in terms of the Schur type functions
associated with the central difference.

\begin{remark}
%%We note some preceding results:
We note some previous results:

\smallskip

\noindent
(1)
We can also express $C^{\mathfrak{o}_N}_k(0)$ and $D^{\mathfrak{sp}_N}_k(0)$
in terms of the symmetrized determinant and 
the symmetrized permanent (\cite{W}, \cite{I5}, \cite{I6}).
Capelli type identities 
in terms of these symmetrized determinant and permanent
are also given in \cite{I3} and \cite{I4}.

\smallskip

\noindent
(2)
%%Analogues of the quantum immanants 
%%in $U(\mathfrak{o}_N)$ and $U(\mathfrak{sp}_N)$ are studied in \cite{OO2}.
%%We will see that these can be regarded as a generalization of 
%%$C^{\mathfrak{o}_N}_k(0)$ and $D^{\mathfrak{sp}_N}_k(0)$
%%in Section~\ref{sec:eigenvalues}.
Analogues of the quantum immanants 
in $U(\mathfrak{o}_N)$ and $U(\mathfrak{sp}_N)$ are studied in \cite{OO2}.
We will see in Section~\ref{sec:eigenvalues}
that these can be regarded as a generalization of 
$C^{\mathfrak{o}_N}_k(0)$ and $D^{\mathfrak{sp}_N}_k(0)$.
\end{remark}
%
%%%%%%%%%%%%%%%%%%%%%%%%%%%%%%%%%%%%%%%%%%%%%%%%%%%%%%%%%%%%%%%%%%%%%%%%%%
%
\section{The eigenvalues of Capelli type elements}\label{sec:eigenvalues}
%
%%%%%%%%%%%%%%%%%%%%%%%%%%%%%%%%%%%%%%%%%%%%%%%%%%%%%%%%%%%%%%%%%%%%%%%%%%
%
Finally we consider concrete examples of Schur type functions
associated with some differences.
When $Q$ is the forward difference,
the associated Schur type functions are equal to the factorial Schur functions.
By the shift of variables,
these are transformed into the shifted Schur functions,
and these are useful to express the eigenvalues of some central elements 
of the universal enveloping algebra of the general linear Lie algebra.
When $Q$ is the central difference,
the associated Schur functions are useful to express the eigenvalues of the central elements 
of the universal enveloping algebras of the orthogonal and symplectic Lie algebras
listed in the previous section.

%%%%%%%%%%%%%%%%%%%%%%%%%%%%%%%%%%%%%%%%%%%%%%%%%%%%%%%%%%%%%%%%%%%%%
\subsection{}
Let us consider the case that $Q$ is equal to the forward difference $\Delta^+$.
In this case, $p^{\Delta^+}_n(x)$ and $p^{*\Delta^+}_n(x)$ are expressed as 
$p^{\Delta^+}_n(x) = x^{\underline{n}}$
and $p^{*\Delta^+}_n(x) = (x-1)^{\underline{n}}$.
The corresponding symmetric functions $e^{\Delta^+}_k$ and $h^{\Delta^+}_k$
are explicitly expressed as follows
($e^{*\Delta^+}_k$ and $h^{*\Delta^+}_k$ are also given by considering the shift of variables).
This expression is essentially equivalent with Corollary~11.3 in \cite{OO1}.

\begin{theorem}\label{thm:explicit_description_of_e_and_h_1}
   \sl 
   We have
   \begin{align*}
      & e^{\Delta^+}_k(x_1,\ldots,x_N) \\ 
	  & \qquad = \sum_{1 \leq i_1 < \cdots < i_k \leq N}
	  (x_{i_1} -N+k-1 + i_1) (x_{i_2} -N+k-2 + i_2) \cdots (x_{i_k} -N + i_k), \\
	  & h^{\Delta^+}_k(x_1,\ldots,x_N) \\
	  & \qquad = \sum_{1 \leq i_1 \leq \cdots \leq i_k \leq N}
	  (x_{i_1} -N-k+1 + i_1) (x_{i_2} -N-k+2 + i_2) \cdots (x_{i_k} -N + i_k).
   \end{align*}
\end{theorem}

\noindent
{\it Proof.}
The first relation is obtained from the relation
$$
   e_N(x_1 - u,\ldots,x_N - u) = (x_1 - u) \cdots (x_N - u).
$$
It suffices to apply $Q_u$ repeatedly
and use the Leibnitz rule for the forward difference:
$\Delta^+ (f(x)g(x)) = \Delta^+ f(x) g(x+1) + f(x) \Delta^+ g(x)$.
The second relation is deduced by induction on~$N$ 
(use Theorems~\ref{thm:gen_fn_of_h_1} and \ref{thm:gen_fn_of_h_2}).
\hfill\qed

\medskip

Using these symmetric functions,
we can describe the eigenvalues of $C^{\frak{gl}_N}_k(u)$ and $D^{\frak{gl}_N}_k(u)$
defined in the previous section as follows:

\begin{theorem}\label{thm:eigenvalues_for_gl}\sl
   For the representation $\pi^{\mathfrak{gl}_N}_{\lambda}$
   of $\mathfrak{gl}_N$ determined by the partition $\lambda=(\lambda_1,\ldots,\lambda_N)$,
   we have
   $$
      \pi^{\mathfrak{gl}_N}_{\lambda}(C^{\mathfrak{gl}_N}_k(u)) 
	  = e^{\Delta^+}_k(l_1,\ldots,l_N ; u), \qquad
	  \pi^{\mathfrak{gl}_N}_{\lambda}(D^{\mathfrak{gl}_N}_k(u))
      = h^{\Delta^+}_k(l_1,\ldots,l_N ; u).
   $$
   Here we put $l_i = \lambda_i + N - i$.
\end{theorem}

%%We can deduce this by the procedure outlined in Section~\ref{sec:capelli_type_elements}
%%noting the triangular decomposition of the general linear Lie algebra
%%and the definitions of the column-determinant and the column-permanent.
We can deduce this by the procedure outlined in Section~\ref{sec:capelli_type_elements}
using the triangular decomposition of the general linear Lie algebra
and the definitions of the column-determinant and the column-permanent.

Theorem~\ref{thm:eigenvalues_for_gl} 
is essentially included in the following description of 
the eigenvalues of the quantum immanants 
in terms of the factorial (shifted) Schur functions
due to \cite{OO1}:
$$
   \pi^{\mathfrak{gl}_N}_{\lambda}(\mathbb{S}_{\mu}) = s^{\Delta^+}_{\mu}(l_1,\ldots,l_N).
$$
Indeed, when $u=0$, Theorem~\ref{thm:eigenvalues_for_gl} 
is a special case of this relation as seen from (\ref{eq:relation_between_C_D_and_QI}).
Moreover, the case of $u \ne 0$ can be also deduced from this by considering
the algebra automorphisms $E_{ij} \mapsto E_{ij} + u \delta_{ij}$
and $E_{ij} \mapsto E_{ij} - u \delta_{ij}$
on $U(\mathfrak{gl}_N)$.

%%%%%%%%%%%%%%%%%%%%%%%%%%%%%%%%%%%%%%%%%%%%%%%%%%%%%%%%%%%%%%%%%%%%%%%%
\subsection{}\label{subsec:central_difference}
Next, let us consider the case of the central difference (namely $Q = \Delta^0$).
In this case, $p^{\Delta^0}_n(x)$ and $p^{*\Delta^0}_n(x)$ are expressed as 
$p^{\Delta^0}_n(x) = x \cdot x^{\overline{\underline{n-1}}}$
and $p^{*\Delta^0}_n(x) = x^{\overline{\underline{n}}}$.
Moreover $e^{\Delta^0}_k$ and $h^{*\Delta^0}_k$ are expressed as follows:

\begin{theorem}\sl
   We have
   \begin{align*}
      & e^{\Delta^0}_k(x_1,\ldots,x_N) \\ 
	  & \quad = \sum_{1 \leq i_1 < \cdots < i_k \leq N}
      (x_{i_1} -\tfrac{N}{2}+\tfrac{k}{2}-1 + i_1) 
	  (x_{i_2} -\tfrac{N}{2}+\tfrac{k}{2}-2 + i_2) \cdots 
	  (x_{i_k} -\tfrac{N}{2}-\tfrac{k}{2} + i_k) \\
	  & \quad = \sum_{1 \leq i_1 < \cdots < i_k \leq N}
      (x_{i_1} +\tfrac{N}{2}-\tfrac{k}{2}+1 - i_1) 
	  (x_{i_2} +\tfrac{N}{2}-\tfrac{k}{2}+2 - i_2) \cdots 
	  (x_{i_k} +\tfrac{N}{2}+\tfrac{k}{2} - i_k), \displaybreak[0]\\
	  & h^{*\Delta^0}_k(x_1,\ldots,x_N) \\
	  & \quad = \sum_{1 \leq i_1 \leq \cdots \leq i_k \leq N}
	  (x_{i_1} -\tfrac{N}{2}-\tfrac{k}{2} + i_1) 
	  (x_{i_2} -\tfrac{N}{2}-\tfrac{k}{2}+1 + i_2) \cdots 
	  (x_{i_k} -\tfrac{N}{2}+\tfrac{k}{2}-1 + i_k) \\
	  & \quad = \sum_{1 \leq i_1 \leq \cdots \leq i_k \leq N}
	  (x_{i_1} +\tfrac{N}{2}+\tfrac{k}{2} - i_1) 
	  (x_{i_2} +\tfrac{N}{2}+\tfrac{k}{2}-1 - i_2) \cdots 
	  (x_{i_k} +\tfrac{N}{2}-\tfrac{k}{2}+1 - i_k).
   \end{align*}
\end{theorem}

The proof of the first relation is almost the same 
as that of Theomem~\ref{thm:explicit_description_of_e_and_h_1}.
Here, we use the Leibnitz rule for the central difference
$$
   \Delta^0 (f(x)g(x)) = \Delta^0 f(x) g(x+\tfrac{1}{2}) + f(x-\tfrac{1}{2}) \Delta^0 g(x).
$$
The second relation is deduced from the second formula in 
Theomem~\ref{thm:explicit_description_of_e_and_h_1}
by replacing $x_i$ with $x'_i = x_i + \frac{N}{2}+\frac{k}{2}$.
Indeed, using elementary row operations, we have
$$
   \begin{vmatrix}
      x_1^{\overline{\underline{N-1+k}}} & \ldots & x_N^{\overline{\underline{N-1+k}}} \\
	  x_1^{\overline{\underline{N-2}}} & \ldots & x_N^{\overline{\underline{N-2}}} \\
	  \vdots & & \vdots \\
	  x_1^{\overline{\underline{0}}} & \ldots & x_N^{\overline{\underline{0}}} 
   \end{vmatrix}
   = \begin{vmatrix}
      x_1^{\prime\underline{N-1+k}} & \ldots & x_N^{\prime\underline{N-1+k}} \\
	  x_1^{\prime\underline{N-2}} & \ldots & x_N^{\prime\underline{N-2}} \\
	  \vdots & & \vdots \\
	  x_1^{\prime\underline{0}} & \ldots & x_N^{\prime\underline{0}} 
   \end{vmatrix},
$$
and this means the relation
$$
   \tilde{s}_{(k)}^{*\Delta^0}(x_1,\ldots,x_N)
   = \tilde{s}_{(k)}^{\Delta^+}(x'_1,\ldots,x'_N).
$$

%%The second relation is deduced from the second relation in 
%%Theomem~\ref{thm:explicit_description_of_e_and_h_1}
%%by replacing $x_i$ by $x'_i = x_i + (k+N)/2$:
%%$$
%%   \tilde{s}_{(k)}^{*\Delta^0}(x_1,\ldots,x_N)
%%   = \begin{vmatrix}
%%      x_1^{\overline{\underline{N-1+k}}} & \ldots & x_N^{\overline{\underline{N-1+k}}} \\
%%	  x_1^{\overline{\underline{N-2}}} & \ldots & x_N^{\overline{\underline{N-2}}} \\
%%	  \vdots & & \vdots \\
%%	  x_1^{\overline{\underline{0}}} & \ldots & x_N^{\overline{\underline{0}}} 
%%   \end{vmatrix}
%%   = \begin{vmatrix}
%%      x_1^{\prime\underline{N-1+k}} & \ldots & x_N^{\prime\underline{N-1+k}} \\
%%	  x_1^{\prime\underline{N-2}} & \ldots & x_N^{\prime\underline{N-2}} \\
%%	  \vdots & & \vdots \\
%%	  x_1^{\prime\underline{0}} & \ldots & x_N^{\prime\underline{0}} 
%%   \end{vmatrix}
%%   = \tilde{s}_{(k)}^{\Delta^+}(x'_1,\ldots,x'_N).
%%$$
%%Here the second equality is deduced by elementary row operations.

Using these functions $e_k^{\Delta^0}$ and $h_k^{*\Delta^0}$, 
%%Using these functions, 
we can express the eigenvalues of the Capelli type central elements 
of the universal enveloping algebras of the orthogonal and symplectic Lie algebras.

First, we have the following relation for $C^{\frak{o}_N}_k(u)$:

\begin{theorem}\label{thm:eigenvalues_for_o}\sl
   For the irreducible representation $\pi^{\mathfrak{o}_N}_{\lambda}$ of $\mathfrak{o}_N$
   determined by the partition $\lambda=(\lambda_1,\ldots,\lambda_{[n]})$,
   we have
   $$
      \pi^{\mathfrak{o}_N}_{\lambda}(C^{\mathfrak{o}_N}_k(u)) 
	  = e^{\Delta^0}_k(l_1,\ldots,l_N ; u).
   $$
   Here we define $l_1,\ldots,l_N$ as follows.
   First, for $1 \leq i \leq n$, we put $l_i = \lambda_i + n - i$
%%   {\rm (}namely we consider the $\rho$-shift{\rm )}.
   {\rm (}that is, we consider the $\rho$-shift{\rm )}.
   Next, we put $l_{n+1} = -l_{n}, \ldots, l_{N} = -l_1$ when $N$ is even,
   and we put $l_{n^{\dagger}} = 0$,
   $l_{n^{\dagger}+1} = -l_{n^{\dagger}-1}, \ldots, l_{N} = -l_1$ 
   with $n^{\dagger} = \frac{N+1}{2}$ when $N$ is odd.
\end{theorem}

Similarly we have the following relation for $D^{\frak{sp}_N}_k(u)$:

\begin{theorem}\label{thm:eigenvalues_for_sp}\sl
   For the irreducible representation $\pi^{\mathfrak{sp}_N}_{\lambda}$ of $\mathfrak{sp}_N$
   determined by the partition $\lambda=(\lambda_1,\ldots,\lambda_n)$,
   we have 
   $$
      \pi^{\mathfrak{sp}_N}_{\lambda}(D^{\mathfrak{sp}_N}_k(u)) 
	  = h^{*\Delta^0}_k(l_1,\ldots,l_N ; u).
   $$
   Here $l_i$ is defined as follows:
   we put $l_i = \lambda_i +n+1-i$ for $1 \leq i \leq n$
%%   {\rm (}namely we consider the $\rho$-shift{\rm )},
   {\rm (}that is, we consider the $\rho$-shift{\rm )},
   and put $l_{n+1} = -l_{n}, \ldots, l_{N} = -l_1$ for $n+1 \leq i \leq N$.
\end{theorem}

\begin{remark}
Factorial powers and differences were key tools
in the study of Capelli type elements,
%%and various relations were given 
and various relations for them have been given 
([I1--6], \cite{IU}, and [U1--5]).
The formulas
in Sections~\ref{sec:expansions_of_e_and_h} and~\ref{sec:generating_functions}
of this article
can be regarded as natural generalizations
of these relations.
\end{remark}

%%%%%%%%%%%%%%%%%%%%%%%%%%%%%%%%%%%%%%%%%%%%%%%%%%%%%%%%%%%%%%%%%%%%%
\subsection{}\label{subsec:relation_with_OO2}
The functions $e^{\Delta^0}_k$ and $h^{* \Delta^0}_k$ 
are related with analogues of the shifted Schur functions
due to Okounkov and Olshanski.
Moreover, from these relations and 
the calculation of the eigenvalues in Theorems~\ref{thm:eigenvalues_for_o} 
and \ref{thm:eigenvalues_for_sp},
we see that $C^{\mathfrak{o}_N}_k(0)$ and $D^{\mathfrak{sp}_N}_k(0)$ 
are equal to 
special cases of analogues of the quantum immanants.
Let us see these relations.

In \cite{OO2}, 
Okounkov and Olshanski introduced 
analogues of the quantum immanants
in the universal enveloping algebras of
the simple Lie algebras of types B, C, and D,
%%and denoted these elements by the symbol $\mathbb{T}_{\mu}$.
and denoted these elements by $\mathbb{T}_{\mu}$.
Let us write these as $\mathbb{T}^B_{\mu}$, $\mathbb{T}^C_{\mu}$, 
and $\mathbb{T}^D_{\mu}$,
when we need to indicate the type of the Lie algebra.
Noting that the eigenvalues of these central elements 
are polynomials in $l_1^2,\ldots,l_n^2$,
Okounkov and Olshanski also defined a class of functions $t^*_{\mu}$ 
corresponding to the types B, C, and D
by the relation 
\begin{align}\label{eq:definition_of_t_star}
   \pi_{\lambda}(\mathbb{T}_{\mu}) = t^*_{\mu}(\lambda_1,\ldots,\lambda_n),
\end{align}
where $\pi_{\lambda}$ is the irreducible representation
of the Lie algebra 
determined by the partition $\lambda$.
%%We can regard this $t^*_{\mu}$ 
We can regard $t^*_{\mu}$ 
as an analogue of the shifted Schur functions,
and this is expressed as follows (Lemma-Definition~2.4 in \cite{OO2}):
\begin{align*}
   t^{*B}_{\mu}(\lambda_1,\ldots,\lambda_n) 
   &= s_{\mu}((\lambda_1+n-\tfrac{1}{2})^2,(\lambda_2+n-\tfrac{3}{2})^2,
   \ldots,(\lambda_n+\tfrac{1}{2})^2 \,|\, 
   (\tfrac{1}{2})^2, (\tfrac{3}{2})^2, (\tfrac{5}{2})^2, \ldots), \\
   t^{*C}_{\mu}(\lambda_1,\ldots,\lambda_n) 
   &= s_{\mu}((\lambda_1+n)^2,(\lambda_2+n-1)^2,\ldots,(\lambda_n+1)^2 
   \,|\, 1^2, 2^2, 3^2, \ldots), \\
   t^{*D}_{\mu}(\lambda_1,\ldots,\lambda_n) 
   &= s_{\mu}((\lambda_1+n-1)^2,(\lambda_2+n-2)^2,\ldots,(\lambda_n+0)^2 
   \,|\, 0^2, 1^2, 2^2,\ldots).
\end{align*}
%%Here, the superscripts $B$, $C$, and $D$
%%mean the type of the corresponding Lie algebra.
Here, the superscripts $B$, $C$, and $D$
indicate the type of the corresponding Lie algebra.
Moreover, the right hand sides are the ``generalized factorial Schur functions'' 
(see \cite{OO2} for the definition).

The connection between these functions and 
$e^{\Delta_0}_k$, $e^{*\Delta_0}_k$, $h^{\Delta_0}_k$, $h^{*\Delta_0}_k$ 
is seen from the following relations:
\begin{align*}
   e^{\Delta^0}_k(l_1,\ldots,l_n,-l_n,\ldots,-l_1) 
   & = s_{(1^k)}(l_1^2,\ldots,l_n^2 \,|\, 0^2,1^2,\ldots), \\
   e^{*\Delta^0}_k(l_1,\ldots,l_n,-l_n,\ldots,-l_1) 
   & = s_{(1^k)}(l_1^2,\ldots,l_n^2 \,|\, (\tfrac{1}{2})^2,(\tfrac{3}{2})^2,\ldots), \\
   e^{*\Delta^0}_k(l_1,\ldots,l_n,0,-l_n,\ldots,-l_1) 
   & = s_{(1^k)}(l_1^2,\ldots,l_n^2 \,|\, 1^2,2^2,\ldots), \\
   h^{*\Delta^0}_k(l_1,\ldots,l_n,-l_n,\ldots,-l_1) 
   & = s_{(k)}(l_1^2,\ldots,l_n^2 \,|\, 1^2,2^2,\ldots), \\
   h^{\Delta^0}_k(l_1,\ldots,l_n,-l_n,\ldots,-l_1) 
   & = s_{(k)}(l_1^2,\ldots,l_n^2 \,|\, (\tfrac{1}{2})^2,(\tfrac{3}{2})^2,\ldots). 
\end{align*}
%%These relations themselves follow from Theorems~\ref{thm:gen_fn_of_e}, 
%%\ref{thm:gen_fn_of_h_1}, and~\ref{thm:gen_fn_of_h_2}.
These relations themselves follow from Theorems~\ref{thm:gen_fn_of_e}--\ref{thm:gen_fn_of_h_2}.
The left hand sides of the second and fourth equalities
are also equal to
$$
   e^{\Delta^0}_k(l_1,\ldots,l_n,0,-l_n,\ldots,-l_1),\qquad
   h^{\Delta^0}_k(l_1,\ldots,l_n,0,-l_n,\ldots,-l_1),
$$
respectively
(recall Proposition~\ref{prop:substitute_0}).

From these relations, 
we can rewrite Theorem~\ref{thm:eigenvalues_for_o} in the case of $u=0$ as
$$
   \pi^{\mathfrak{o}_{2n}}_{\lambda}(C^{\mathfrak{o}_{2n}}_k(0))
   =
   t^{* D}_{(1^k)}(\lambda_1,\ldots,\lambda_n), \qquad
   \pi^{\mathfrak{o}_{2n+1}}_{\lambda}(C^{\mathfrak{o}_{2n+1}}_k(0))
   =
   t^{* B}_{(1^k)}(\lambda_1,\ldots,\lambda_n).
$$
Similarly, 
Theorem~\ref{thm:eigenvalues_for_sp} 
in the case of $u=0$ can be rewritten as
$$
   \pi^{\mathfrak{sp}_{2n}}_{\lambda}(D^{\mathfrak{sp}_{2n}}_k(0))
   =
   t^{* C}_{(k)}(\lambda_1,\ldots,\lambda_n).
$$
Combining these with (\ref{eq:definition_of_t_star}),
we have the following theorem:

\begin{theorem}\label{thm:relation_with_analogues_of_quantum_immanants}\sl
   We have
   $$
      C^{\mathfrak{o}_{2n}}_k(0) = \mathbb{T}^D_{(1^k)}, \qquad
      C^{\mathfrak{o}_{2n+1}}_k(0) = \mathbb{T}^B_{(1^k)}, \qquad
      D^{\mathfrak{sp}_N}_k(0) = \mathbb{T}^C_{(k)}.
   $$
\end{theorem}

Indeed, the eigenvalues of both sides are equal.
 
This theorem means that
we can express $\mathbb{T}^D_{(1^k)}$,
$\mathbb{T}^B_{(1^k)}$,
and $\mathbb{T}^C_{(k)}$ in terms of the column-determinant
and the column-permanent.
However, an explicit descripsion of $\mathbb{T}_{\mu}$ in terms of a certain noncommutative
determinant type function is not given for general $\mu$.
Thus 
the three elements in Theorem~\ref{thm:relation_with_analogues_of_quantum_immanants}
are lucky exceptions.
In \cite{MN}, these 
three elements are studied in more detail,
and Capelli type identities for the dual pair $(O_M,Sp_N)$ are given.

In Theorems~\ref{thm:eigenvalues_for_o} and \ref{thm:eigenvalues_for_sp},
we expressed the eigenvalues of $C^{\mathfrak{o}_N}_k(u)$ and $D^{\mathfrak{sp}_N}_k(u)$
naturally in terms of the functions $e^{\Delta^0}_k$ and $h^{*\Delta^0}_k$ 
for general $u$.
This is an advantage of these functions over the function $t^*_{\mu}$.
%%Indeed, we can not describe these eigenvalues  
Indeed, we cannot describe these eigenvalues  
in terms of $t^{*B}_{\mu}$, $t^{*C}_{\mu}$, and $t^{*D}_{\mu}$ 
so simply.
%%However it should be noted that, even if $u \ne 0$,
However, it should be noted that, even if $u \ne 0$,
we can express 
the eigenvalue of $C^{\mathfrak{o}_N}_k(u)$
as a linear combination of 
$t^{*D}_{(1^0)},t^{*D}_{(1^1)},\ldots,t^{*D}_{(1^k)}$
or $t^{*B}_{(1^0)},t^{*B}_{(1^1)},\ldots,t^{*B}_{(1^k)}$
using Theorem~\ref{thm:expansion_of_e}.
Similarly, 
the eigenvalue of $D^{\mathfrak{sp}_N}_k(u)$
can be expressed 
as a linear combination of 
$t^{*C}_{(0)},t^{*C}_{(1)},\ldots,t^{*C}_{(k)}$
by using Theorem~\ref{thm:expansion_of_h_1}.

In the various relations in this article,
there was a mysterious duality in the exchanges $s \leftrightarrow s^*$ and 
$\lambda \leftrightarrow \lambda'$.
It is also mysterious that 
the functions $e$ and $h^*$ played more important roles than $e^*$ and $h$
(note that we can rewrite Theorem~\ref{thm:eigenvalues_for_gl}
in terms of $h^*$ replacing $u$ by $u-1$).
This puzzling phenomenon seems to be related to the following fact:
central elements in $U(\mathfrak{o}_N)$ (respectively, $U(\mathfrak{sp}_N)$) 
expressed in terms of the column-permanent (respectively, the column-determinant)
are not known.
The author hopes that 
%%the theoretical background of these phenomena will be transparent,
the theoretical background of these phenomena will become transparent,
and the Schur type functions in this article will be useful to study
the analogues of the quantum immanants in $U(\mathfrak{o}_N)$ and $U(\mathfrak{sp}_N)$
(especially to give their explicit description).
%
%%%%%%%%%%%%%%%%%%%%%%%%%%%%%%%%%%%%%%%%%%%%%%%%%%%%%%%%%%%%%%%%%%%%%%%%%%
%
% References
%
%%%%%%%%%%%%%%%%%%%%%%%%%%%%%%%%%%%%%%%%%%%%%%%%%%%%%%%%%%%%%%%%%%%%%%%%%%
%

\end{document}